\documentclass[12pt,a4paper]{article}
\usepackage{amsmath}
\usepackage{amssymb}
\usepackage{amsthm}
\usepackage{amsfonts}
\usepackage{latexsym}

\theoremstyle{plain}
\newtheorem{thm}{Theorem}[section]
\newtheorem{cor}{Corollary}[section]
\newtheorem{lem}{Lemma}[section]
\newtheorem{prop}{Proposition}[section]
\theoremstyle{remark}

\theoremstyle{definition}
\newtheorem{defn}{Definition}[section]

\newtheorem{rem}{Remark}[section]

\title{Equivariant asymptotics for Toeplitz operators}
\author{Roberto Paoletti\footnote{\noindent{\bf Address:}
Dipartimento di Matematica e Applicazioni, Universit\`a degli Studi
di Milano Bicocca, Via R. Cozzi 53, 20125 Milano,
Italy; {\bf e-mail}: roberto.paoletti@unimib.it }}
\date{}

\begin{document}

\maketitle

\begin{abstract}
In recent years, the Tian-Zelditch asymptotic expansion for the equivariant components of the Szeg\"{o}
kernel of a polarized complex projective manifold, and its subsequent generalizations
in terms of scaling limits, have played an important role in algebraic, symplectic,
and differential geometry. A natural question is whether there exist generalizations in which the
projector onto the spaces of holomorphic sections can be replaced by
the projector onto more general (non-complete) linear series. One case that lends itself to such analysis,
and which is natural from the point of view of geometric quantization,
is given by the linear series determined by imposing spectral bounds on an invariant self-adjoint Toeplitz operator.
In this paper we focus on the asymptotics of the spectral projectors associated to slowly
shrinking spectral bands.

\end{abstract}

\section{Introduction}

Let $M$ be a d-dimensional complex projective manifold, $A$ an ample line bundle
on it. Suppose that $h$ is an Hermitian metric on $A$, and that the unique connection
compatible with the Hermitian and holomorphic structures has normalized
curvature $\Theta=-2i\,\omega$, where $\omega$ is a K\"{a}hler form on $M$.
Then $dV_M=:\big(1/\mathrm{d}!\big)\,\omega^{\wedge \mathrm{d}}$ is a volume form
on $M$, with total volume $\mathrm{vol}(M)=:\left(\pi^\mathrm{d}/\mathrm{d}!\right)\,
\int_M c_1(A)^\mathrm{d}$.

The dual line bundle $A^{-1}=A^\vee$ naturally inherits an Hermitian structure, and
the unit circle bundle $X\subset A^\vee$ is a principal $S^1$-bundle
on $M$; let $\pi:X\rightarrow M$ denote the projection.
Then $A$ is the line bundle associated to $X$ and the standard representation of
$S^1$ on $\mathbb{C}$. In particular, for every $k\in \mathbb{Z}$ there are natural isomorphisms
$\mathcal{C}^\infty\left(M,A^{\otimes k}\right)\cong \mathcal{C}^\infty (X)_k$, where the left
hand side is the space of smooth global sections of $A^{\otimes k}$, and the right hand side
denotes the $k$-th isotype in $\mathcal{C}^\infty (X)$ for the $S^1$-action.

The normalized connection form, $\alpha\in \Omega^1(X)$, is a contact structure on $X$,
hence $d\mu_X=:(1/2\pi)\, \alpha\wedge \pi^*(dV_M)$ is a volume form on $X$. With these
choices, the above isomorphisms are unitary with respect to the natural Hermitian structures.
For $k=0,1,2,\ldots$ the space $H^0\left(M,A^{\otimes k}\right)$
of global holomorphic sections of $A^{\otimes k}$
corresponds to the $k$-th isotype $H(X)_k=H(X)\cap \mathcal{C}^\infty (X)_k$
of the Hardy space $H(X)\subset L^2(X)$; with this in mind,
we shall occasionally implicitly identify $H^0\left(M,A^{\otimes k}\right)$
and $H(X)_k$.

The \textit{Szeg\"{o} projector} is the orthogonal projector $\Pi:L^2(X)\rightarrow H(X)$; $\Pi$ extends
to a linear operator $\mathcal{D}'(X)\rightarrow H(X)$ (we shall implicitly identify
functions, densities and half-densities by the given choices).
The object of this paper are certain asymptotic properties of $S^1$-invariant
Toeplitz operators on $X$, that is, operators
of the form $T=\Pi\circ P\circ \Pi$, where $P$ is an $S^1$-invariant pseudodifferential
operator of classical type.
The $S^1$-invariance of $P$ implies that $T$ preserves the decomposition into
$S^1$-isotypes, and the asymptotics in point refer to the Fourier decomposition. More precisely, let
$\Pi_k:L^2(X)\rightarrow H(X)_k$
be the orthogonal projector; then
$T=\bigoplus _kT_k$, where $T_k=\Pi_k\circ P\circ \Pi_k$,
and we are interested in the asymptotics of the spectral function of
$T_k$ along the diagonal of $X$, for $k\rightarrow +\infty$.

In algebro-geometric terms, we shall thus study the local asymptotics of
families of (possibly non-complete)
linear series
determined by spectral bounds imposed by an invariant Toeplitz operator.
This way of defining a linear series is unorthodox
in algebraic geometry, but seems quite natural from the perspective of geometric quantization,
where it corresponds to imposing un upper bound on, say, the total energy of the system.

Hence, on the one hand this work specializes the local study of Toeplitz operators of \cite{p}
to the equivariant context; in this sense, the main Theorem below is
an equivariant version of
the local Weyl law for Toeplitz operator of \cite{p}.
On the other hand, it may be seen as a generalization to the Toeplitz context of the
Tian-Zelditch asymptotic expansion \cite{t}, \cite{z}, and of the scaling limits in \cite{bsz}, \cite{sz},
where the full Szeg\"{o} kernel $\Pi$ is replaced by
the spectral function of $T$.

The theme of this paper is related as well to the theory of \cite{bpu}, which also
deals with the asymptotics of certain spectral projectors associated to Toeplitz operators.
The focus in \cite{bpu} is on vector subspaces
associated to fast narrowing bands of energy levels, and on the asymptotics at fixed pairs of points;
in particular, it is proved in \cite{bpu}
that at a point $x\in X$ at energy level $E$ eigensections of energy differing by at most
$c\,k^{-1}$ from $E$ give a contribution to the full Szeg\"{o} kernel which grows 
like an appropriate multiple of $k^{\mathrm{d}-1/2}$ (while the full Szeg\"{o} kernel
grows like $k^{\mathrm{d}}$).
This raises the natural question to determine how the contribution of wider energy bands
relates to the full Szeg\"{o} kernel at $x$.
Thus we consider energy bands that shrink at a relatively slow rate
(for example, fixed bands), and
estimate their contribution to the full Szeg\"{o} kernel at
pairs of points converging to each other
at a controlled rate as $k\rightarrow +\infty$.
In this discussion, \lq energy\rq\,
is thought of as a Toeplitz operator of order zero, while the result of the paper
will be phrased in terms of first order operators; one passes from one to the
other by composition with the elliptic Toeplitz operator associated with
the circle action on $X$, which turns a fixed spectral energy band into one expanding at a rate
linear with $k$.

To describe our results, it is order to recall some notation from \cite{bg}.

\begin{defn}
\label{defn:toeplitz-machinery}
Let $X$ and $\Pi$ be as above.
\begin{itemize}
  \item A Toeplitz operator of order $m\in \mathbb{Z}$ on $X$ is an
  operator $T:\mathcal{D}'(X)\rightarrow \mathcal{D}'(X)$ of the form $T=\Pi\circ P\circ \Pi$,
  where $P$ is a pseudodifferential operator of classical type of order $m$ on $X$.
  \item Let
  $$
  \Sigma=:\Big\{(x,r\,\alpha_x):\,x\in X,\,r>0\Big\}\subset T^*X\setminus\{0\}.
  $$
  If $T:\mathcal{D}'(X)\rightarrow \mathcal{D}'(X)$ is a Toeplitz operator, its symbol
  $\sigma_T :\Sigma\rightarrow \mathbb{C}$ is the restriction of the symbol of $P$. Thus
  $\sigma_T$ is real if $T$ is self-adjoint.
  \item The \textit{reduced symbol} $\varsigma_T\in \mathcal{C}^\infty(X)$ is
  $\varsigma_T(x)=:\sigma_T\big(x,\alpha_x\big)$ ($x\in X$). If $T$ is $S^1$-invariant,
  $\varsigma_T$ may be regarded as a smooth function on $M$.
\end{itemize}
\end{defn}

Suppose then that $T$ is a  first order $S^1$-invariant self-adjoint  Toeplitz operator;
set $a_T=:\min \varsigma_T$, $A_T=:\max \varsigma_T$.
For every $k$, we may regard $T_k$ in a natural manner as a self-adjoint endomorphism
$T_k:H(X)_k\rightarrow H(X)_k$. With this interpretation,
let $\lambda_{k1}\le \cdots \le \lambda_{kN_k}$ be the
eigenvalues of $T_k$, repeated according to multiplicity.
Then for every $j$
\begin{equation}
\label{eqn:asympt-spectral-bd}
a_T\,k+O(1)\le \lambda_{kj}\le A_T\,k+O(1)
\end{equation}
as $k\rightarrow +\infty$ (a proof will be given below).
For every $k=0,1,2,\ldots$, we can find an orthonormal basis $\big(e_{kj}\big)$ of $H(X)_k$
such that $e_{kj}$ is an eigenvector of $T_k$ with eigenvalue $\lambda_{kj}$,
for every $j=1,\ldots,N_k$.

\begin{defn}
\label{defn:equiv-spectral-fctn}
The level-k spectral function of $T$ is
$$
\mathcal{T}_k\Big(\lambda,x',x''\Big)=:\sum _{j:\,\lambda_{kj}\le \lambda}e_{kj}\left(x'\right)\,\overline{e_{kj}\left(x''\right)}
\,\,\,\,\,\,\,(\lambda\in \mathbb{R},\,x',x''\in X).
$$
\end{defn}

Thus $\mathcal{T}_k\big(\lambda,\cdot,\cdot)\in \mathcal{C}^\infty(X\times X)$
is the kernel of the orthogonal projector onto the span $V_\lambda^{(k)}\subseteq H(X)_k$ of the eigenspaces
of $T_k$ corresponding to eigenvalues $\le \lambda$; in particular, it does not depend on $\big(e_{kj}\big)$.
The asymptotic bound (\ref{eqn:asympt-spectral-bd}) motivates restricting attention to the
asymptotics of $\mathcal{T}_k\big(\lambda\,k,\cdot,\cdot)$ as $k\rightarrow +\infty$, where $\lambda$ is fixed.

Our result will be expressed in Heisenberg local coordinates centered at a given $x\in X$ \cite{sz}.
This implies choosing first a system of preferred local coordinates on $M$ centered at $m=:\pi(x)$,
meaning that the symplectic and complex structures on the tangent space $T_mM$ are the standard
ones, and then a preferred local frame $e_L$ of $A$ centered at $m$,
meaning that the \lq Hessian\,\rq $\nabla^2e_L$ at $m$ is as expected in the local Heisenberg model (see \cite{sz}
for a precise discussion). By \cite{sz}, in Heisenberg local coordinates the scaling limits of Szeg\"{o}
kernels exhibit a universal nature, and the point of this work is that in certain ranges
the same holds of the equivariant
spectral functions of invariant Toeplitz operators.

For any $\ell\in \mathbb{N}$ and $\delta>0$, let
$B_\ell(\mathbf{0},\delta)\subseteq \mathbb{R}^\ell$ be the open ball
of radius $\delta$ centered at the origin.
Following \cite{sz}, if $x\in X$ and $\mathfrak{h}:(-\pi,\pi)\times B_{2\mathrm{d}}(\mathbf{0},\delta)\rightarrow X$
is a system of Heisenberg local coordinates centered at $x$, we shall
set $x+(\theta,\mathbf{v})=:\mathfrak{h}(\theta,\mathbf{v})$, and occasionally
$x+\mathbf{v}=\mathfrak{h}(0,\mathbf{v})$. If $r_\vartheta:X\rightarrow X$ is the action of
$e^{i\vartheta}\in S^1$, then $r_\vartheta\big(x+(\theta,\mathbf{v})\big)=
x+(\vartheta+\theta,\mathbf{v})$. The given system of preferred local coordinates determines
a unitary isomorphism $\mathbb{C}^\mathrm{d}\cong T_mM$, and with this understanding this
notation will be applied to suitably small $\mathbf{v},\mathbf{w}\in T_mM$.

Finally, we need a further piece of notation from \cite{sz}.

\begin{defn}
Let $H$ be the Hermitian structure on $M$ determined by $\omega$; thus, $\omega=-\Im(H)$.
Let $\|\cdot\|$ be the norm associated to $H$.
If $m\in M$ and $\mathbf{w},\mathbf{v}\in T_mM$, we shall let
\begin{eqnarray*}
\psi_2(\mathbf{w},\mathbf{v})&=:&i\,\Im\big(H_m(\mathbf{w},\mathbf{v})\big)-\frac 12\,\|\mathbf{w}-\mathbf{v}\|_m^2\\
&=&H_m(\mathbf{w},\mathbf{v})-\frac 12\,\left(\|\mathbf{w}\|_m^2+\|\mathbf{v}\|_m^2\right).
\end{eqnarray*}
\end{defn}

We can now state:

\begin{thm}
\label{thm:main}
Let $T$ be a first order $S^1$-invariant self-adjoint Toeplitz operator on $X$.
Suppose $x\in X$, $m=:\pi(x)$. Suppose $0\le \xi<1/2$, $c>0$,
and $\varpi\ge 0$ satisfies $\varpi\le 1/6$ and $\varpi< 1/2-\xi$.
Let $e_k\in \mathbb{R}$
be a sequence such that $e_k>c\,k^{-\xi}$, $k=1,2,\ldots$.
\begin{enumerate}
  \item Uniformly in  $\lambda\le \varsigma_T(m)-e_k$,
  in $\mathbf{w},\mathbf{v}\in T_mM$
  with $\max\{\|\mathbf{w}\|,\,\|\mathbf{v}\|\}\lesssim k^{\varpi}$
  and in $\theta,\theta'\in (-\pi,\pi)$, as $k\rightarrow+\infty$ we have
  $$
  \mathcal{T}_k\left(\lambda\,k,x+
  \left(\theta,\frac{\mathbf{w}}{\sqrt{k}}\right),x+\left(\theta',\frac{\mathbf{v}}{\sqrt{k}}\right)\right)=
  O\left(k^{-\infty}\right).
  $$
  \item Uniformly in  $\lambda\ge \varsigma_T(m)+e_k$,
  in $\mathbf{w},\mathbf{v}\in T_mM$
  with $\max\{\|\mathbf{w}\|,\,\|\mathbf{v}\|\}\lesssim k^{\varpi}$,
  and in $\theta,\theta'\in (-\pi,\pi)$, as $k\rightarrow+\infty$ for every $N=1,2,\ldots$
  we have
  \begin{eqnarray}
  \label{eqn:spectral-szego}
  \lefteqn{\mathcal{T}_k\left(\lambda\,k,x+
  \left(\theta,\frac{\mathbf{w}}{\sqrt{k}}\right),x+\left(\theta',\frac{\mathbf{v}}{\sqrt{k}}\right)\right)}\\
  &=& \Pi_k\left(x+
  \left(\theta,\frac{\mathbf{w}}{\sqrt{k}}\right),
  x+\left(\theta',\frac{\mathbf{v}}{\sqrt{k}}\right)\right)+O\left(k^{-\infty}\right).\nonumber
  \end{eqnarray}
\end{enumerate}
\end{thm}

The proof will combine classical arguments in the study of spectral functions of pseudodifferential
operators \cite{h}, \cite{gs} and microlocal tecnhiques from \cite{z}, \cite{bsz}, \cite{sz}
revolving around the description of the Szeg\"{o} kernel as a Fourier integral \cite{bs};
some basic results about Toeplitz operators from \cite{bg} will be key ingredients in the proof.

The following remarks are in order.

First, by replacing $T$ with $-T$, in Theorem \ref{thm:main} statement 1 about lower bands
turns into a statement about upper bands. Thus, statement 2 is really an expansion regarding any slowly
shrinking
intermediate energy band containing $\varsigma_T(m)$ in its interior. In fact, we may use an orthonormal
basis of eigenvectors of $T$ to estimate the asymptotics of the full equivariant Szeg\"{o} kernel
$\Pi_k$ in \cite{sz}; thus the latter may be written
as the sum of three terms, one from a lower band, one from an intermediate band containing $\varsigma_T(m)$ in its interior,
and one from an upper band, as above. Hence, 2 is a consequence of 1.

To state the previous point explicitly, for $\lambda_1<\lambda_2$ define
\begin{eqnarray*}
\mathcal{T}_k\left(\lambda_1,\lambda_2;x',x''\right)=:
\mathcal{T}_k\left(\lambda_2,x',x''\right)-\mathcal{T}_k\left(\lambda_1,x',x''\right).
\end{eqnarray*}
That is, $\mathcal{T}_k\left(\lambda_1,\lambda_2;x',x''\right)$ is the kernel of the
orthogonal projector onto the span of the eigenspaces corresponding to eigenvalues comprised
in the half-open band $(\lambda_1,\lambda_2]$.

\begin{cor}
In the hypothesis of Theorem \ref{thm:main}, as $k\rightarrow +\infty$ we have
\begin{eqnarray*}
\lefteqn{\mathcal{T}_k\left(k\big(\varsigma_T(m)-e_k\big),k\big(\varsigma_T(m)+e_k\big);x+
  \left(\theta,\frac{\mathbf{w}}{\sqrt{k}}\right),
  x+\left(\theta',\frac{\mathbf{v}}{\sqrt{k}}\right)\right)}\\
&=&\Pi_k\left(x+
  \left(\theta,\frac{\mathbf{w}}{\sqrt{k}}\right),
  x+\left(\theta',\frac{\mathbf{v}}{\sqrt{k}}\right)\right)+O\left(k^{-\infty}\right).
  \end{eqnarray*}
\end{cor}

Secondly, by \cite{sz} the scaling limit of the Szeg\"{o} kernel on the right hand side
of (\ref{eqn:spectral-szego}) has
a \lq large ball\rq\, asymptotic expansion for $k\rightarrow +\infty$, with a universal leading term
$$
e^{ik(\theta-\theta')+\psi_2(\mathbf{w},\mathbf{v})}\,
\left(\frac k\pi\right)^{\mathrm{d}}.
$$
The same then holds of the spectral function on the left hand side of (\ref{eqn:spectral-szego}).
A succinct direct derivation
using stationary phase techniques is given below.

One motivation for this work is to extend the study of the asymptotic properties of complete linear
series to linear series defined
by spectral bounds on Toeplitz operators. Let us give a couple of immediate applications.

Under mild assumptions, the rate of growth of the dimension of the linear series $V_{\lambda k}^{(k)}$
is governed by the volume of the locus of phase space where the Hamiltonian $\varsigma_T\le \lambda$.
More precisely,
let $M_{< \lambda}=:\big\{m\in M:\varsigma_T(m)< \lambda\big\}$.

\begin{cor}
\label{cor:asympt-dim}
Assume that $\lambda$ is a regular value of $\varsigma_T$.
Then
$$
\lim_{k\rightarrow +\infty}\left(\frac{\pi}{k}\right)^{\mathrm{d}}
\,\dim \left(V_{\lambda k}^{(k)}\right)=
\mathrm{vol}\big(M_{< \lambda}\big).
$$
\end{cor}
\noindent
A similar statement obviously holds for intermediate bands.

Next, we consider the asymptotics of the rational maps associated to $V_{\lambda k}^{(k)}$.
If $L$ is a line bundle on $M$ and $V\subseteq H^0\left(M,L\right)$ is a vector space of holomorphic sections,
the \textit{base locus}
$\mathrm{Bs}(V)$ of the linear series $|V|$ is the common zero locus of all sections in $V$. Thus $m\not\in \mathrm{Bs}(V)$
if and only if there exists $s\in V$ such that $s(m)\neq 0$.

The following Corollary is proved as in the study of the full linear series in \cite{z}
(actually establishing asymptotic isometry on compact subsets of $M_{< \lambda}$).

\begin{cor}
\label{cor:base-point}
If $\varsigma_T(x)<\lambda$ then $\pi(x)\not\in \mathrm{Bs}\left(V_{\lambda k}^{(k)}\right)$,
and the rational map induced by the linear series $\left|V_{\lambda k}^{(k)}\right|$ is immersive at
$\pi(x)$, for all $k\ge k_x$.
\end{cor}


\section{Preliminaries}

In this section we shall quickly put things in perspective by proving the asymptotic estimate
(\ref{eqn:asympt-spectral-bd}). Let $\langle \cdot,\cdot\rangle$
be the Hermitian product on $L^2(X)$.

If $f\in \mathcal{C}^\infty(X)$, we shall denote multiplication by
$f$ by
$M_f:\mathcal{D}'(X)\rightarrow \mathcal{D}'(X)$, $g\mapsto f\,g$. If $f\in \mathcal{C}^    \infty(M)$,
we shall regard it in the natural manner as an $S^1$-invariant function on $X$, and denote by
$T_f=\Pi\circ M_f\circ \Pi$ the associated invariant zero order
Toeplitz operator, and by $T_f^{(k)}:H(X)_k\rightarrow H(X)_k$
the endomorphisms induced by restriction.
If $f\in \mathcal{C}^    \infty(M)$ is real, every $T_f^{(k)}$ is self-adjoint; let
$\lambda_{k1}\le \ldots\le \lambda_{kN_k}$ be its eigenvalues, repeated according to multiplicity.

\begin{lem}
\label{lem:multiplication-case}
Suppose $f\in \mathcal{C}^\infty(M)$ is real and let $a_f=:\min f$, $A_f=:\max f$. Then
$a_f\le \lambda_{kj}\le A_f$, for $k=0,1,2,\ldots$ and $1\le j\le N_k$.
\end{lem}

\textit{Proof.} If $\sigma\in L^2(X)$ has unit norm, then
\begin{equation*}
a_f\le \int_X f\,\sigma\,\overline{\sigma}\,d\mu_X\le A_f.
\end{equation*}
On the other hand, because $\Pi$ is self-adjoint if in addition
$\sigma\in H(X)$ we have
\begin{eqnarray*}
\langle T_f(\sigma),\sigma\rangle =\langle M_f(\sigma),\sigma\rangle=
\int_X f\,\sigma\,\overline{\sigma}\,d\mu_X.
\end{eqnarray*}
The statement follows.

\hfill Q.E.D.

\bigskip

Let $\partial _\theta$ be the vector field on $X$ generating
the $S^1$-action, and set $D=:-i\,\partial _\theta$.
If $f\in \mathcal{C}^\infty(M)$, $\widetilde{T}_f=:\Pi\circ (D\circ M_f)\circ \Pi=
D\circ T_f$ is an invariant first order Toeplitz operator, self-adjoint if
$f$ is real. By Lemma \ref{lem:multiplication-case},
its eigenvalues $\widetilde{\lambda}_{kj}=k\,\lambda_{kj}$
satisfy $a_f\,k\le \widetilde{\lambda}_{kj}\le k\,A_f$.

Since $T_D=:\Pi\circ D\circ \Pi$ is an elliptic invariant Toeplitz operator of degree $1$,
there exists an invariant Toeplitz operator of degree $-1$ such that $E\circ D=\Pi+S$,
where $S$ is smoothing and invariant; in particular, the norm of $S$ on $H(X)_k$ is
$O\left(k^{-\infty}\right)$. Since $D$ on $H(X)_k$ is
$k\,\mathrm{id}$, $E$ induces endomorphisms $E_k:H(X)_k\rightarrow H(X)_k$
satisfying $k\,E_k=\mathrm{id}+O\left(k^{-\infty}\right)$, hence
$E_k=k^{-1}\,\mathrm{id}+O\left(k^{-\infty}\right)$.

\begin{lem}
\label{lem:asympt-general-case-zero}
Let $T$ be a zero order $S^1$-invariant self-adjoint Toeplitz operator,
$f=:\varsigma _T$. Then $a_f+O\left(k^{-1}\right)\le\lambda_{kj}\le
A_f+O\left(k^{-1}\right)$ as $k\rightarrow +\infty$.
\end{lem}

\textit{Proof.}
We have $T=T_f+R$, where $R$ is an invariant self-adjoint Toeplitz operator of degree
$-1$. Now $R=\Pi\circ R=E\circ (D\circ R)+R'$, where $R'=-S\circ R$ is smoothing.
Since $D\circ R$ is a Toeplitz operator of degree $0$, it is bounded in norm, and
therefore the previous discussion implies that $R$ is $O\left(k^{-1}\right)$ on $H(X)_k$.
The statement follows from this and Lemma \ref{lem:multiplication-case}

\hfill Q.E.D.

\begin{cor}
\label{cor::asympt-general-case-zero}
Let $T$ be a first order $S^1$-invariant self-adjoint Toeplitz operator,
$f=:\varsigma _T$. Then $a_f\,k+O\left(1\right)\le\lambda_{kj}\le
A_f\,k+O\left(1\right)$ as $k\rightarrow +\infty$.
\end{cor}


\section{Proof of Theorem \ref{thm:main}}

As noted in the introduction, we need only prove 1.

To begin with, we may reduce the proof to the case $\varsigma_T\ge 1$.
For if $C\ge 0$ is such that $C+\varsigma_T\ge 1$, and
$T'=:T+C\,D$, then
$\varsigma_{T'}=\varsigma_T+C$, and an orthonormal
basis $(e_{kj})$ of eigenvectors of $T$ with eigenvalues
$\lambda_{kj}$ is also a basis of eigenvectors of $T'$,
with eigenvalues $\lambda_{kj}'=:\lambda_{kj}+C\,k$.
Hence the spectral functions $\mathcal{T}$ and $\mathcal{T}'$
of $T$ and $T'$ are related by
$\mathcal{T}'_k\big((\lambda+C)\,k,x',x''\big)=\mathcal{T}_k\big(\lambda\,k,
x',x''\big)$, and $\lambda \lessgtr \varsigma_T(x)\,\Leftrightarrow\,
\lambda+C\lessgtr \varsigma_{T'}(x)$.
Thus the asymptotic expansion for $T'$ implies the one for $T$.

Furthermore, by construction the equivariant spectral function $\mathcal{T}_k$
only involves the $k$-th isotype for the $S^1$-action, therefore
\begin{eqnarray*}
\mathcal{T}_k\left(\lambda\,k,x+\left(\theta,\frac{\mathbf{w}}{\sqrt{k}}\right),
x+\left(\theta',\frac{\mathbf{v}}{\sqrt{k}}\right)\right)=
e^{ik(\theta-\theta')}\,\mathcal{T}_k\left(\lambda\,k,x+\frac{\mathbf{w}}{\sqrt{k}},
x+\frac{\mathbf{v}}{\sqrt{k}}\right).
\end{eqnarray*}
Hence we may assume without loss that $\theta=\theta'=0$.

As a further reduction, it suffices to prove the theorem when $\mathbf{w}=\mathbf{0}$.
To see this, recall that preferred and Heisenberg local coordinates may be
deformed smoothly with the reference point.
More precisely, there exist first an open neighborhood $U\subseteq M$ of $m$ and for
every $m'\in U$ preferred local coordinates $\mathfrak{p}_{m'}$ centered at, and smoothly depending
on, $m'$; next, there exists for every $x'\in \pi^{-1}(U)$ a system of Heisenberg local
coordinates $\mathfrak{h}_{x'}$ centered at and smoothly depending
on $x'$, and such that the system of preferred local coordinates underlying $\mathfrak{h}_{x'}$
is $\mathfrak{p} _{\pi\left(x'\right)}$. If $\delta>0$ is sufficiently small and
$\mathbf{w}\in B_{2\mathrm{d}}(\mathbf{0};\delta)$, we shall let
$m'+\mathbf{w}=:\mathfrak{p}_{m'}(\mathbf{w})$ if $m'\in U$,
and $x'+\mathbf{w}=:\mathfrak{h}_{x'}(0,\mathbf{w})$
 if $x'\in \pi^{-1}(U)$.

If $\mathbf{w},\,\mathbf{v}\in B_{2\mathrm{d}}(\mathbf{0};\delta)$ then
$(m+\mathbf{w})+\mathbf{v}=m+A(\mathbf{w},\mathbf{v})$ for a certain $\mathbb{C}^{\mathrm{d}}$-valued
smooth function $A$;
by Taylor expansion, for $N=1,2,\ldots$ we get
\begin{equation}
\label{eqn:change-of-prefer-coord}
(m+\mathbf{w})+\mathbf{v}=m+\left(\sum_{j=1}^N R_j(\mathbf{w},\mathbf{v})
+O\left(\max\big\{\|\mathbf{w}\|,\|\mathbf{v}\|\big\}^{N+1}\right)\right),
\end{equation}
where $R_j$ is a homogeneous $\mathbb{C}^{\mathrm{d}}$-valued polynomial of degree $j$.
Actually, $R_1(\mathbf{w},\mathbf{v})=\mathbf{w}+\mathbf{v}$ in (\ref{eqn:change-of-prefer-coord}).
More precisely, working in rescaled coordinates for future reference we have

\begin{lem}
\label{lem:compare-sum-M}
Suppose $\mathbf{w},\mathbf{v}\in B_{2\mathrm{d}}(\mathbf{0},R)$. Then as $k\rightarrow +\infty$
$$
\left(m+\frac{\mathbf{w}}{\sqrt{k}}\right)+\frac{\mathbf{v}}{\sqrt{k}}=
m+\left(\frac{1}{\sqrt{k}}\,(\mathbf{w}+\mathbf{v})+O\left(\frac{R^2}{k}\right)\right).
$$
\end{lem}

\textit{Proof.} For suitably small $\upsilon>0$ and $m'\in U$, let
$\mathfrak{P}_{m'}=:\mathfrak{p}_m^{-1}\circ \mathfrak{p}_{m'}:
B_{2\mathrm{d}}(\mathbf{0};\upsilon)\rightarrow \mathbb{R}^{2\mathrm{d}}$.
Thus
$d_\mathbf{0}\mathfrak{P}_{m'}=\mathrm{id}_{\mathbb{R}^{2\mathrm{d}}}
+O\left(\|m'\|\right)$, where $\|m'\|=:\left\|\mathfrak{p}_m^{-1}(m')\right\|$.
Therefore,
$$
d_\mathbf{0}\mathfrak{P}_{m+\frac{\mathbf{w}}{\sqrt{k}}}
\left(\frac{\mathbf{v}}{\sqrt{k}}\right)=\frac{\mathbf{v}}{\sqrt{k}}+
O\left(\frac{R^2}{k}\right).
$$
By construction,
$$
\mathfrak{p}_m^{-1}\left(m+\frac{1}{\sqrt{k}}\,(\mathbf{w}+\mathbf{v})\right)=
\frac{1}{\sqrt{k}}\,(\mathbf{w}+\mathbf{v}).
$$
On the other hand, again by construction,
$$
\mathfrak{P}_{m+\frac{\mathbf{w}}{\sqrt{k}}}(\mathbf{0})=\mathfrak{p}_m^{-1}\left(
\mathfrak{p}_{m+\frac{\mathbf{w}}{\sqrt{k}}}(\mathbf{0})\right)=
\mathfrak{p}_m^{-1}\left(m+\frac{\mathbf{w}}{\sqrt{k}}\right)=\frac{\mathbf{w}}{\sqrt{k}}.
$$
 Thus,
\begin{eqnarray*}
\lefteqn{\mathfrak{p}_m^{-1}\left(
\left(m+\frac{\mathbf{w}}{\sqrt{k}}\right)+\frac{\mathbf{v}}{\sqrt{k}}\right)=
\mathfrak{p}_m^{-1}\circ \mathfrak{p}_{m+\frac{\mathbf{w}}{\sqrt{k}}}
\left(\frac{\mathbf{v}}{\sqrt{k}}\right)}\\
&=&\mathfrak{P}_{m+\frac{\mathbf{w}}{\sqrt{k}}}\left(\frac{\mathbf{v}}{\sqrt{k}}\right)
=\mathfrak{P}_{m+\frac{\mathbf{w}}{\sqrt{k}}}(\mathbf{0})+d_\mathbf{0}\mathfrak{P}_{m+\frac{\mathbf{w}}{\sqrt{k}}}
\left(\frac{\mathbf{v}}{\sqrt{k}}\right)+O\left(\frac{R^2}{k}\right)\\
&=&\frac{\mathbf{w}}{\sqrt{k}}+\frac{\mathbf{v}}{\sqrt{k}}+O\left(\frac{R^2}{k}\right).
\end{eqnarray*}
\hfill Q.E.D.

\bigskip

We now lift this comparison to Heisenberg coordinates on $X$.
If $\mathbf{w},\,\mathbf{v}\in B_{2\mathrm{d}}(\mathbf{0};\delta)$ then
$m+\mathbf{v}=(m+\mathbf{w})+\big(\theta(\mathbf{w},\mathbf{v}),B(\mathbf{w},\mathbf{v})\big)$,
for suitable smooth real and $\mathbb{C}^{\mathrm{d}}$-valued smooth functions $\theta$ and $B$,
respectively. Taylor expansion then yields
\begin{eqnarray}
\theta(\mathbf{w},\mathbf{v})&=&\sum_{j=1}^N \theta_j(\mathbf{w},\mathbf{v})+
O\left(\max\big\{\|\mathbf{w}\|,\|\mathbf{v}\|\big\}^{N+1}\right),\\
B(\mathbf{w},\mathbf{v})&=&\sum_{j=1}^N B_j(\mathbf{w},\mathbf{v})
+O\left(\max\big\{\|\mathbf{w}\|,\|\mathbf{v}\|\big\}^{N+1}\right);
\label{eqn:change-of-heisen-coord}
\end{eqnarray}
here $\theta_j$ and $B_j$ are homogenous of degree $j$.
In fact, $B_1(\mathbf{w},\mathbf{v})=\mathbf{v}-\mathbf{w}$,
$\theta_1=0$, and $\theta_2 (\mathbf{w},\mathbf{v})=\omega_0(\mathbf{\mathbf{w}},\mathbf{v})$.
Since it is not essential in the following argument, we state this without proof (in rescaled
coordinates):

\begin{lem}
\label{lem:compare-sum-X}
In the hypothesis of Lemma \ref{lem:compare-sum-M}, as $k\rightarrow +\infty$ we have
$$
x+\frac{\mathbf{v}}{\sqrt{k}}=\left(x+\frac{\mathbf{w}}{\sqrt{k}}\right)+
\left(\theta\left(\frac{1}{\sqrt{k}}\right),\frac{1}{\sqrt{k}}\,(\mathbf{v}-\mathbf{w})+O\left(\frac{R^2}{k}\right)\right),
$$
where $\theta(\upsilon)=\upsilon^2\,\omega_0(\mathbf{\mathbf{w}},\mathbf{v})+O\left(R^3\,\upsilon^3\right)$
as $\upsilon\rightarrow 0$.
\end{lem}

Here $\omega_0$ is the standard symplectic structure on $\mathbb{R}^{2\mathrm{d}}\cong
\mathbb{C}^{2\mathrm{d}}$.
If $\mathbf{w},\,\mathbf{v}$ are interpreted as elements in
$T_mM$, then $\omega_0(\mathbf{\mathbf{w}},\mathbf{v})$ should be replaced by
$\omega_m(\mathbf{\mathbf{w}},\mathbf{v})$.

Suppose now that the statement of Theorem \ref{thm:main} has been proved when
$\mathbf{w}=\mathbf{0}$.
If $\|\mathbf{w}\|,\,\|\mathbf{v}\|\lesssim k^{\varpi}$, let us set
$m(k)=:m+\mathbf{w}/\sqrt{k}$.
Then $\big|\varsigma_T(m)-\varsigma_T\big(m(k)\big)\big| \lesssim k^{\varpi-1/2}$, thus
if $\lambda< \varsigma_T(m)-e_k$ then $\lambda<\varsigma_T\big(m(k)\big)-
e_k/2$ for $k\gg 0$, since by assumption $-\xi>\varpi-1/2$.
In this range, using $x+\mathbf{w}/\sqrt{k}$ as reference point we obtain:
\begin{eqnarray*}
\lefteqn{\mathcal{T}_k\left(\lambda\,k,x+\frac{\mathbf{w}}{\sqrt{k}},
x+\frac{\mathbf{v}}{\sqrt{k}}\right)}\\
&=&\mathcal{T}_k\left(\lambda\,k,x+\frac{\mathbf{w}}{\sqrt{k}},
\left(x+\frac{\mathbf{w}}{\sqrt{k}}\right)
+\left(
\theta\left(\frac{1}{\sqrt{k}}\right),
\left(\frac{1}{\sqrt{k}}\,(\mathbf{v}-\mathbf{w})+\frac 1k\,\rho\left(\frac 1{\sqrt{k}}\right)\right)
\right)
\right)\nonumber\\
&=&e^{-ik\theta(1/\sqrt{k})}\,\mathcal{T}_k\left(\lambda\,k,x+\frac{\mathbf{w}}{\sqrt{k}},
\left(x+\frac{\mathbf{w}}{\sqrt{k}}\right)
+\frac{1}{\sqrt{k}}\,(\mathbf{v}-\mathbf{w})+\frac 1k\,\rho\left(\frac 1{\sqrt{k}}\right)\right)\\
&=&O\left(k^{-\infty}\right).
\end{eqnarray*}

We now prove the Theorem assuming $\varsigma_T\ge 1$,
$\theta=\theta'=0$, $\mathbf{w}=\mathbf{0}$.

In view of Lemma 12.1 of \cite{bg}, perhaps averaging we can find a first order $S^1$-invariant
self-adjoint pseudodifferential operator $Q$ on $X$, such that $T=\Pi\circ Q\circ \Pi$,
$[Q,\Pi]=0$, and with positive principal symbol $q:T^*X\setminus \{0\}\rightarrow (0,+\infty)$.
In Heisenberg local coordinates, the $S^1$-action is a translation in
$\theta$, therefore by $S^1$-invariance
$q\big(x+(\vartheta+\theta,\mathbf{v})\big)=  q\big(x+(\theta,\mathbf{v})\big)$.

Since $q>0$, we have $Q\ge -c\,\mathrm{id}$ for some $c\in \mathbb{R}$. Thus
$Q'=:Q+(c+1)\,\mathrm{id}\ge \mathrm{id}$. Let $T'=:\Pi\circ Q'\circ \Pi$; then
$T'=T+(c+1)\,\Pi$, and $\varsigma_T=\varsigma_{T'}$.
Now the orthogonal basis $(e_{kj})$ of eigenvectors of $T$, with eigenvalues
$\lambda_{kj}$, is also an orthogonal basis of eigenvectors of $T'$, with eigenvalues
$\lambda_{kj}'=\lambda_{kj}+(c+1)$. Hence, the spectral functions $\mathcal{T}$
and $\mathcal{T}'$ of $T$ and $T'$ are
related by the equality
$$
\mathcal{T}\left(\lambda\,k,x',x''\right)=\mathcal{T}'\left(\left(\lambda+\frac 1k\,
(c+1)\right)\cdot k,x',x''\right).
$$
But if $\lambda\le \varsigma_T(m)- e_k$,
then $\lambda+
(c+1)/k<\varsigma_T(m)- e_k/2$ for all $k\gg 0$. Hence if the statement holds for
$T$ if it holds for $T'$. Thus we are reduced to proving the theorem under the further
assumption $Q\ge \mathrm{id}$.

Given a sufficiently small
$\epsilon>0$ and $\chi\in \mathcal{C}^\infty_0\big((-\epsilon,\epsilon)\big)$, for any suitable family of
operators $A(\tau)$, after \cite{gs} we set
$$
A_\chi=:\int _{-\epsilon}^\epsilon \chi(\tau)\,A(\tau)\,d\tau.
$$
This may be applied to the 1-parameter group of unitary operators $U(\tau)=:e^{i\tau Q}$,
and to their equivariant Toeplitz contrations $S^{(k)}(\tau)=:\Pi_k\circ U(\tau)\circ \Pi_k$; furthermore, we shall
replace $\chi$ with $\chi\,e^{-ik\lambda (\cdot)}$. Then
$S_{\chi\,e^{-ik\lambda (\cdot)}}^{(k)}=U_{\chi\,e^{-ik\lambda (\cdot)}}\circ \Pi_k$ is a smoothing operator, with
Schwartz kernel
\begin{equation}
\label{eqn:S-chi-k}
S_{\chi\,e^{-ik\lambda (\cdot)}}^{(k)}\left(x',x''\right)=
\sum_{j=1}^{N_k}\widehat{\chi}\big(\lambda \,k-\lambda_{kj}\big)\,e_{kj}\left(x'\right)\cdot \overline{e_{kj}\left(x''\right)}
\end{equation}
(similar functions have been studied in \cite{bpu}). Now in this construction we replace $\chi$ with
its rescaling $\chi_k(\tau)=:\chi\left(k^\xi\,\tau\right)$. Then the Fourier transform of
$\chi_k$ is $\widehat{\chi}_k(s)=(1/k^\xi)\,\widehat{\chi}(s/k^\xi)$ ($s\in \mathbb{R}$), hence
\begin{equation}
\label{eqn:S-chi-k-k}
S_{\chi_k\,e^{-ik\lambda (\cdot)}}^{(k)}\left(x',x''\right)=
\frac{1}{k^\xi}\sum_{j=1}^{N_k}\widehat{\chi}\left(\frac{1}{k^{\xi}}\,\big(\lambda \,k-\lambda_{kj}\big)\right)
\,e_{kj}\left(x'\right)\cdot \overline{e_{kj}\left(x''\right)}
\end{equation}

We introduce complex measures on the real line
$$
\mu_{\mathcal{T}^{(k)}}=\mu_{\mathcal{T}^{(k)}}^{(x,\mathbf{v})}=:
\sum_{j=1}^{N_k}e_{kj}\left(x\right)\cdot \overline{
e_{kj}\left(x+\frac{\mathbf{v}}{\sqrt{k}}\right)}\,
\delta_{\lambda_{kj}},
$$
where $\delta_t$ is the delta function at $t$. Thus,
$\mathcal{T}^{(k)}\left(\lambda',x,x+\mathbf{v}/\sqrt{k}\right)=
\int_{-\infty}^{\lambda'} d\mu_{\mathcal{T}^{(k)}}$, $\forall \,\lambda'\in \mathbb{R}$.
Then
\begin{equation}
\label{eqn:S-chi-k-as-integral}
S_{\chi_k\,e^{-ik\lambda (\cdot)}}^{(k)}\left(x,x+\frac{\mathbf{v}}{\sqrt{k}}\right)=
\int_{-\infty}^{+\infty}\widehat{\chi}_k\big(\lambda \,k-\eta\big)\,d\mu_{\mathcal{T}^{(k)}}(\eta).
\end{equation}

Define $G_k(\eta)=:\int _{-\infty}^\eta \widehat{\chi}_k(b)\,db$ ($\eta\in \mathbb{R}$).
Then $\int _{-\infty}^{+\infty}G_k(k\lambda-\eta)\, d\mu_{\mathcal{T}^{(k)}}(\eta)$
may be computed in two different manners, and comparing the results will
yield the stated asymptotic expansion.

Let us embark on the first computation. We have:
\begin{eqnarray}
\label{eqn:comput-1st-1}
\lefteqn{\int _{-\infty}^{+\infty}G_k(k\lambda-\eta)\, d\mu_{\mathcal{T}^{(k)}}(\eta)}\\
&=&\int _{-\infty}^{+\infty}\left[\int_{-\infty}^{k\lambda-\eta}\widehat{\chi}_k(b)\,db\right]\,
d\mu_{\mathcal{T}^{(k)}}(\eta)=
\int _{-\infty}^{+\infty}\left[\int_{-\infty}^{k\lambda}\widehat{\chi}_k(b-\eta)\,db\right]\,
d\mu_{\mathcal{T}^{(k)}}(\eta).\nonumber
\end{eqnarray}

\begin{lem}
\label{lem:commute-integrals}
We have
$$
\int _{-\infty}^{+\infty}\left[\int_{-\infty}^{k\lambda}\widehat{\chi}_k(b-\eta)\,db\right]\,
d\mu_{\mathcal{T}^{(k)}}(\eta)=
\int_{-\infty}^{k\lambda}\left[\int _{-\infty}^{+\infty}
\widehat{\chi}_k(b-\eta)\,d\mu_{\mathcal{T}^{(k)}}(\eta)\right]
\,db.
$$
\end{lem}

\textit{Proof.} Leaving dependence on $x$ and $\mathbf{v}$ implicit, define
$\beta_{kj}\in (-\pi,\pi]$ and $r_{kj}>0$ by the equality
$
r_{kj}\,e^{i\beta_{kj}}=e_{kj}\left(x\right)\cdot \overline{
e_{kj}\left(x+\mathbf{v}/\sqrt{k}\right)}$.
The total variation of $\mu_{\mathcal{T}^{(k)}}$ is then
$\left|\mu_{\mathcal{T}^{(k)}}\right|=\sum_j r_{kj}\,\delta_{\lambda_{kj}}$.
Let $\beta_k\in \mathcal{C}^\infty(\mathbb{R})$ be a real function such that
$\beta_k\left(\lambda_{kj}\right)=\beta_{kj}$ for every $j$, and set $h_k=e^{i\beta_k}$.
Thus $\mu_{\mathcal{T}^{(k)}}=h_k\,\left|\mu_{\mathcal{T}^{(k)}}\right|$, and
\begin{eqnarray*}
\lefteqn{\int _{-\infty}^{+\infty}\left[\int_{-\infty}^{k\lambda}\widehat{\chi}_k(b-\eta)\,db\right]\,
d\mu_{\mathcal{T}^{(k)}}(\eta)}\\
&=&
\int _{-\infty}^{+\infty}\left[\int_{-\infty}^{k\lambda}\widehat{\chi}_k(b-\eta)\,h_k(\eta)\,db\right]\,
d\left|\mu_{\mathcal{T}^{(k)}}\right|(\eta).
\end{eqnarray*}
Let $C=:\left\|\widehat{\chi}\right\|_{L^1}=\left\|\widehat{\chi}_k\right\|_{L^1}$; then
\begin{eqnarray*}
\lefteqn{\int _{-\infty}^{+\infty}\left[\int_{-\infty}^{k\lambda}\left|\widehat{\chi}_k(b-\eta)\,h_k(\eta)\right|\,db\right]\,
d\left|\mu_{\mathcal{T}^{(k)}}\right|(\eta)}\\
&=&\int _{-\infty}^{+\infty}\left[\int_{-\infty}^{k\lambda}\left|\widehat{\chi}_k(b-\eta)\right|\,db\right]\,
d\left|\mu_{\mathcal{T}^{(k)}}\right|(\eta)\le
C\,\int _{-\infty}^{+\infty}d\left|\mu_{\mathcal{T}^{(k)}}\right|(\eta)\\
&=&C\,\sum_{j=1}^{N_k}r_{kj}<+\infty.
\end{eqnarray*}
The Fubini-Tonelli Theorem then implies
\begin{eqnarray*}
\lefteqn{\int _{-\infty}^{+\infty}\left[\int_{-\infty}^{k\lambda}\widehat{\chi}_k(b-\eta)\,db\right]\,
d\mu_{\mathcal{T}^{(k)}}(\eta)}\\
&=&\int _{-\infty}^{+\infty}\left[\int_{-\infty}^{k\lambda}\widehat{\chi}_k(b-\eta)\,h_k(\eta)\,db\right]\,
d\left|\mu_{\mathcal{T}^{(k)}}\right|(\eta)\\
&=&\int _{-\infty}^{k\lambda}
\left[\int_{-\infty}^{+\infty}\widehat{\chi}_k(b-\eta)\,h_k(\eta)\,
d\left|\mu_{\mathcal{T}^{(k)}}\right|(\eta)\right]\,db\\
&=&\int _{-\infty}^{k\lambda}
\left[\int_{-\infty}^{+\infty}\widehat{\chi}_k(b-\eta)\,
d\mu_{\mathcal{T}^{(k)}}(\eta)\right]\,db.
\end{eqnarray*}

\hfill Q.E.D.

\bigskip

Applying Lemma \ref{lem:commute-integrals} and performing the change of
variable $b\rightsquigarrow k\,b$, we get from (\ref{eqn:comput-1st-1})
\begin{eqnarray}
\label{eqn:1-st-comput-key-reduction}
\lefteqn{\int _{-\infty}^{+\infty}G_k(k\lambda-\eta)\, d\mu_{\mathcal{T}^{(k)}}(\eta)=
\int_{-\infty}^{k\lambda}\left[\int _{-\infty}^{+\infty}
\widehat{\chi}_k(b-\eta)\,d\mu_{\mathcal{T}^{(k)}}(\eta)\right]
\,db}\\
&=&k\,\int_{-\infty}^{\lambda}\left[\int _{-\infty}^{+\infty}
\widehat{\chi}_k\big(k\,b-\eta\big)\,d\mu_{\mathcal{T}^{(k)}}(\eta)\right]
\,db =k\,\int_{-\infty}^{\lambda}S^{(k)}_{\chi_k\,e^{-ikb(\cdot)}}
\left(x,x+\frac{\mathbf{v}}{\sqrt{k}}\right)\,db.
\nonumber
\end{eqnarray}
We shall now consider the asymptotics of the latter integral.
Let $a_{\varsigma_T}=:\min \varsigma_T$,
$A_{\varsigma_T}=:\max \varsigma_T$; thus $a_{\varsigma_T}\ge 1$.

\begin{lem}
\label{lem:irrelevant-interval}
For every $N=1,2,\ldots$ there exists a constant $C_N>0$ such that
$$
\left|S^{(k)}_{\chi_k\,e^{-ikb(\cdot)}}
\left(x',x''\right)\right|\le C_N\,k^{2\mathrm{d}-N(1-\xi)-\xi}\,(|b|+1)^{-N}
$$
for every $\left(x',x''\right)\in X\times X$ and $b\not\in (1/2,A_{\varsigma_T}+1)$.
\end{lem}

Before commencing the proof, we notice that since $\chi$ is compactly
supported, $\widehat{\chi}$ is of rapid decay; therefore, for every $N>0$
there exists $C_N>0$ such that
$\left |\widehat{\chi}(b)\right|<C_N\,\big(1+|b|\big)^{-N}$ for $b\in \mathbb{R}$.
Hence for every $N=1,2,\ldots$ we have for $k\rightarrow +\infty$:
\begin{eqnarray}
\label{eqn:rescaled-rapid-decay}
\lefteqn{\left |\widehat{\chi}_k(b)\right|=\frac{1}{k^\xi}\,\left |\widehat{\chi}\left(\frac{b}{k^\xi}\right)\right|
\le C_{2N}\,\frac{k^{(2N-1)\xi}}{\left(k^\xi+|b|\right)^{2N}}}\\
&\le&C_{2N}\,\frac{k^{(2N-1)\xi}}{\left[k^{2N\xi}+{2N\choose N}
k^{N\xi}|b|^N\right]}\le C_{2N}'\,\frac{k^{(N-1)\xi}}{1+
|b|^N}\le D_N\,\frac{k^{(N-1)\xi}}{\left(1+
|b|\right)^N}. \nonumber
\end{eqnarray}

\textit{Proof.} By assumption $\varsigma_T\ge 1$, whence
$\lambda_{kj}\ge k+O(1)$ by
Corollary \ref{cor::asympt-general-case-zero}.
If $-1/2\le b\le 1/2$, then $\left|\lambda _{kj}-kb\right|\ge k/3$ for all
$k\gg 0$ and $j=1,\ldots,N_k$. By (\ref{eqn:rescaled-rapid-decay}), for
every $N=1,2,\ldots$ we then have
$$
\left|\widehat{\chi}_k(kb-\lambda_{kj})\right|\le B_N\,k^{-N(1-\xi)-\xi}
\le B_N'\,k^{-N(1-\xi)-\xi}\,(|b|+1)^{-N},
$$ for a constant
$B_N$ independent of $k$ and $j$, and $B_N'=2^N\,B_N$.
If instead $b\not\in (-1/2,A_{\varsigma_T}+1)$,
$\left|\lambda _{kj}-kb\right|\ge C\,k\,(|b|+1)$ for all $k\gg0$ and
some constant $C>0$ independent of $k$ and $j$. Therefore, again by (\ref{eqn:S-chi-k})
we obtain
$$
\left|\widehat{\chi}_k(kb-\lambda_{kj})\right|\le
B_N'\,k^{-N(1-\xi)-\xi}\,(|b|+1)^{-N}
$$
in this range also.

On the other hand, by the Tian-Zelditch asymptotic expansion
$\left|\sigma_{kj}(y)\right|=O\left(k^{\mathrm{d}/2}\right)$ for all $y\in X$,
and by the Riemann-Roch Theorem
$\dim H(X)_k=O\left(k^{\mathrm{d}}\right)$. Therefore, in view of (\ref{eqn:S-chi-k}) we have
$$
\left|S^{(k)}_{\chi_k\,e^{-ikb(\cdot)}}
\left(x',x''\right)\right|\le D_N\,k^{2\mathrm{d}-N(1-\xi)-\xi}\,(|b|+1)^{-N}.
$$

\hfill Q.E.D.

\begin{cor}
\label{cor:irrelevant-interval}
Uniformly in $\left(x',x''\right)\in X\times X$, as $k\rightarrow +\infty$ we have
$$
\int_{-\infty}^{1/2}S^{(k)}_{\chi_k\,e^{-ikb(\cdot)}}
\left(x',x''\right)\,db=O\left(k^{-\infty}\right).
$$
\end{cor}

\begin{rem}
\label{rem:bound-on-lambda}
The same argument implies
$$
\int_{-\infty}^{a_{\varsigma_T}-\delta}S^{(k)}_{\chi_k\,e^{-ikb(\cdot)}}
\left(x',x''\right)\,db=O\left(k^{-\infty}\right)
$$
for any $\delta>0$. Similarly,
$$
\int_{A_{\varsigma_T}+\delta}^{+\infty}S^{(k)}_{\chi_k\,e^{-ikb(\cdot)}}
\left(x',x''\right)\,db=O\left(k^{-\infty}\right).
$$
Henceforth we assume
$
A_f+1>\lambda>1/2$.
\end{rem}

Letting $\sim$ denote equal asymptotics,
(\ref{eqn:1-st-comput-key-reduction}) and Lemma \ref{lem:irrelevant-interval}
imply
\begin{equation}
\label{eqn:1-st-comput-key-asymptotics}
\int _{-\infty}^{+\infty}G_k(k\lambda-\eta)\, d\mu_{\mathcal{T}^{(k)}}(\eta)
\sim k\,\int_{1/2}^{\lambda}S^{(k)}_{\chi_k\,e^{-ikb(\cdot)}}
\left(x,x+\frac{\mathbf{v}}{\sqrt{k}}\right)\,db.
\end{equation}
By definition,
\begin{eqnarray}
\label{eqn:fourier-component-of-S}
\lefteqn{S^{(k)}_{\chi_k\,e^{-ikb(\cdot)}}
\left(x',x''\right)=\frac{1}{2\pi}\,\int_{-\pi}^{\pi}e^{-ik\vartheta}\,S_{\chi_k\,e^{-ikb(\cdot)}}
\left(r_\vartheta\left(x'\right),x''\right)\,d\vartheta}\\
&=&\frac{1}{2\pi}\,\int_{-\pi}^{\pi}\int_{-\epsilon}^\epsilon
e^{-ik(b\tau+\vartheta)}\,\chi_k(\tau)\,\big(U(\tau)\circ \Pi)\left(r_\vartheta\left(x'\right),x''\right)\,d\vartheta\,
d\tau.\nonumber
\end{eqnarray}

As explained in \S 12 of \cite{gs}, if $\epsilon$ is sufficiently small for all
$\tau\in (-\epsilon,\epsilon)$ we can write
$U(\tau)=V(\tau)+R(\tau)$, where  $V(\tau)$ and $R(\tau)$ are as follows.

$V(\tau)$ is a Fourier integral operator, locally of the form
\begin{equation}
\label{eqn:V-tau}
V(\tau)\left(x',x''\right)=\frac{1}{(2\pi)^{2\mathrm{d}+1}}\,\int_{\mathbb{R}^{2\mathrm{d}+1}}
e^{i\left[\varphi\left(\tau,x',\eta\right)-x''\cdot \eta\right]}\,
a\left(\tau,x',x'',\eta\right)\,d\eta;
\end{equation}
here $a(\tau,\cdot,\cdot)\in S^0_{\mathrm{cl}}$, and
\begin{equation}\label{eqn:phase-v-tau}
    \varphi\left(\tau,x',\eta\right)=x'\cdot \eta+\tau\,q\left(x',\eta\right)+
    O\left(\tau^2\right)\cdot \|\eta\|,
\end{equation}
where $x'\cdot \eta$ is the standard Euclidean pairing between the local coordinates of $x'$
and $\eta\in \mathbb{R}^{2\mathrm{d}+1}$.
As discussed in \cite{p},
since the density bundle is trivialized by $d\mu_X$
the initial condition $U(0)=\mathrm{id}$ implies $a\left(0,x',x'',\eta\right)=1/\mathcal{V}\left(x''\right)$,
where $d\mu_X(y)=\mathcal{V}(y)\,dy$ in local coordinates.

$R(\tau)$ is a smooth family of smoothing operators
on $\mathcal{D}'(X)$, parametrized by $\tau$. More precisely, its kernel
$\left(\tau,x',x''\right)\mapsto R\left(\tau,x',x''\right)$ is in
$\mathcal{C}^\infty\big((-\epsilon,\epsilon)\times X\times X\big)$.

Now define $S^{(k)'}_{\chi_k\,e^{-ikb(\cdot)}}$ and $S^{(k)''}_{\chi_k\,e^{-ikb(\cdot)}}$
as in (\ref{eqn:fourier-component-of-S}),
with $U(\tau)$ replaced by $V(\tau)$ and $R(\tau)$ respectively.
Since $U(\tau)=V(\tau)+R(\tau)$, (\ref{eqn:fourier-component-of-S}) implies
$S^{(k)}_{\chi_k\,e^{-ikb(\cdot)}}=S^{(k)'}_{\chi_k\,e^{-ikb(\cdot)}}+S^{(k)''}_{\chi_k\,e^{-ikb(\cdot)}}$.

\begin{lem}
\label{lem:S''-rapid-decrease}
For every $N=1,2,\ldots$, there exist constants $C_N>0$ such that
for all $\left(x',x''\right)\in X\times X$ and $b\in \mathbb{R}$ we have
$$
\left|S^{(k)''}_{\chi_k\,e^{-ikb(\cdot)}}
\left(x',x''\right)\right|
\le C_N\,k^{-N}\,(1+|b|)^{-N}.
$$
\end{lem}

\textit{Proof.}
Since $R(\tau)$ is a smooth family of smoothing operators, the same holds
of $R(\tau)\circ \Pi$. Therefore, the kernel of the latter family defines a smooth function $\mathcal{R}$
on $(-\epsilon,\epsilon)\times X\times X$.

Thus,
$
\widehat{\mathcal{R}}_k\left(\tau,x',x''\right)=:(1/2\pi)\,\int _{-\pi}^{\pi}
e^{-ik\vartheta}\,\mathcal{R}\left(\tau,r_\vartheta\left(x'\right),x''\right)\,d\vartheta
=O\left(k^{-\infty}\right)$ in $\mathcal{C}^j$-norm, uniformly
on $(-\epsilon,\epsilon)\times X\times X$.

By definition,
\begin{equation}
\label{defn:S-ii-fourier-comp}
S^{(k)''}_{\chi_k\,e^{-ikb(\cdot)}}
\left(x',x''\right)=
\int_{-\epsilon}^\epsilon
e^{-ik b\tau}\,\chi_k(\tau)\,\widehat{\mathcal{R}}_k\left(\tau,r_\vartheta\left(x'\right),x''\right)\,
d\tau,
\end{equation}
hence $\left|S^{(k)''}_{\chi_k\,e^{-ikb(\cdot)}}
\left(x',x''\right)\right|\le C_Nk^{-N}$ uniformly in $\left(x',x''\right)\in X\times X$ and $b\in \mathbb{R}$.
Therefore the statement holds for $|b|\le 1$ with $C_N$ replaced by $2^N\,C_N$.

On the other hand, if $|b|\ge 1$ integrating by parts in $d\tau$ in (\ref{defn:S-ii-fourier-comp})
yields
$$
S^{(k)''}_{\chi_k\,e^{-ikb(\cdot)}}
\left(x',x''\right)=\big(i/kb\big)^r\,\int_{-\epsilon}^\epsilon
e^{-ikb\tau}
\big(d/d\tau\big)^r\,\left(\chi_k(\tau)\cdot \widehat{R}_k\right)\left(\tau,x',x''\right)\,d\tau
$$
for $r=1,2,\ldots$. Hence, $\left|S^{(k)''}_{\chi_k\,e^{-ikb(\cdot)}}
\left(x',x''\right)\right|\le |b|^{-r}\,O\left(k^{-\infty}\right)$ for $|b|\ge 1$.

\hfill Q.E.D.

\begin{cor}
\label{lem:S''-rapid-decrease}
Uniformly in $\left(x',x''\right)\in X\times X$ and $-\infty\le\lambda_1<\lambda_2<+\infty$,
as $k\rightarrow +\infty$ we have
$$
\int_{\lambda_1}^{\lambda_2} S^{(k)''}_{\chi_k\,e^{-ikb(\cdot)}}
\left(x',x''\right)\,db=O\left(k^{-\infty}\right).
$$
\end{cor}

Let $P_k:L^2(X)\rightarrow L^2(X)_k$
be the orthogonal projector, given by
$$
P_k(f)(x)=:\frac{1}{2\pi}\, \int_{-\pi}^\pi
e^{-ik\vartheta}\,f\big(r_\vartheta(x)\big)\,d\vartheta,
$$
and
set $V^{(k)}_{\chi_k\,e^{-ikb(\cdot)}}=:P_k\circ V_{\chi_k\,e^{-ikb(\cdot)}}$.
By $S^1$-invariance, $V(\tau)\circ P_k=P_k\circ V(\tau)$, therefore
\begin{eqnarray*}
\lefteqn{
V^{(k)}_{\chi_k\,e^{-ikb(\cdot)}}\circ \Pi}\\
&=&\left(P_k\circ V_{\chi_k\,e^{-ikb(\cdot)}}\right)
\circ \Pi=V_{\chi_k\,e^{-ikb(\cdot)}}\circ \big(P_k\circ \Pi)=V_{\chi_k\,e^{-ikb(\cdot)}}\circ
\Pi_k.
\end{eqnarray*}

On the upshot, in view of (\ref{eqn:1-st-comput-key-asymptotics}),
(\ref{eqn:fourier-component-of-S}), and Lemma \ref{lem:S''-rapid-decrease},
\begin{eqnarray}
\label{eqn:V-tau-only-matters}
\lefteqn{\int _{-\infty}^{+\infty}G_k(k\lambda-\eta)\, d\mu_{\mathcal{T}^{(k)}}(\eta)
\sim k\,\int_{1/2}^{\lambda}S^{(k)'}_{\chi_k\,e^{-ikb(\cdot)}}
\left(x,x+\frac{\mathbf{v}}{\sqrt{k}}\right)\,db}\\
&=&k\,\int_{1/2}^{\lambda}
\left(V^{(k)}_{\chi_k\,e^{-ikb(\cdot)}}\circ \Pi\right)
\left(x,x+\frac{\mathbf{v}}{\sqrt{k}}\right)\,db
=k\,\int_{1/2}^{\lambda}
\left(V_{\chi_k\,e^{-ikb(\cdot)}}\circ \Pi_k\right)
\left(x,x+\frac{\mathbf{v}}{\sqrt{k}}\right)\,db.                                 \nonumber
\end{eqnarray}

If $\beta>0$, let $B_M(m,\beta)\subseteq M$ be the $\beta$-neighborhood of $m$
in the Riemannian distance $\mathrm{dist}_M$.
If $\beta$ is small enough, $M_1=:B_M(m,2\beta)\subseteq \mathfrak{p}\big(
B_{2\mathrm{d}}(\mathbf{0},\delta)\big)$, where $\mathfrak{p}:
B_{2\mathrm{d}}(\mathbf{0},\delta)\rightarrow M$ is the given preferred coordinate chart
centered at $m$.
Set $M_2=:B_M(m,\beta)^c$. Then $\mathcal{M}=\{M_1,M_2\}$ is an open cover of $M$;
let $\{\varrho^{(l)}\}_{l=1}^2$ be a smooth partition of unity subordinate to $\mathcal{M}$.
Using $\mathfrak{p}$, we can pull-back $\varrho^{(l)}$
to $\mathbb{C}^\mathrm{d}\cong \mathbb{R}^{2\mathrm{d}}$, $\mathbf{v}\mapsto \varrho^{(j)}(m+\mathbf{v})$
(recall that $m+\mathbf{v}=\mathfrak{p}(\mathbf{v})$).
Since $\mathfrak{p}$ is a local isometry at the origin, we may assume
$\mathrm{supp}\left(\varrho^{(1)}\circ \mathfrak{p}\right)\subseteq B_{2\mathrm{d}}(\mathbf{0},3\beta)$,
$\mathrm{supp}\left(\varrho^{(2)}\circ \mathfrak{p}\right)\subseteq B_{2\mathrm{d}}(\mathbf{0},\beta/2)^c$.

For $k=1,2,\ldots$ and $l=1,2$, define $\varrho^{(lk)}:M\rightarrow \mathbb{R}$ by setting
\begin{equation}
\label{eqn:def-of-varrho-lk}
\varrho^{(lk)}(m+\mathbf{\mathbf{u}})=:\varrho^{(l)}\left(m+k^{1/2-\varpi}\,\frac \beta 6
\,\mathbf{u}\right).
\end{equation}
Then $\left\{\varrho^{(lk)}\right\}_{l=1}^2$ is a smooth partition of unity on
$M$ for each $k$, and
$$
\mathrm{supp}\left(\varrho^{(1k)}\circ \mathfrak{p}\right)\subseteq B_{2\mathrm{d}}\left(\mathbf{0},18\,k^{\varpi-1/2}\right),
\,\,\,\,\,\,
\mathrm{supp}\left(\varrho^{(2k)}\circ \mathfrak{p}\right)\subseteq B_{2\mathrm{d}}\left(\mathbf{0},3\,k^{\varpi-1/2}\right)^c.
$$
In particular,
\begin{equation}
\label{eqn:bound-on-distance-rho-2}
m'\in \mathrm{supp}\left(\varrho^{(2k)}\right)
\,\,\,\,\,\,
\Rightarrow \,\,\,\,\,\,
\mathrm{dist}_M\left(m,m'\right)\ge \frac 32\,k^{\varpi-1/2}.
\end{equation}

Now $\left\{\varrho^{(lk)}\circ \pi\right\}_{l=1}^2$ is a smooth $S^1$-invariant
partition of unity on $X$; let us define $\Pi^{(lk)}\in \mathcal{D}'(X\times X)$ by
setting $\Pi^{(lk)}\left(x',x''\right)=:\varrho^{(lk)}\big(\pi\left(x'\right)\big)\,
\Pi\left(x',x''\right)$. Then $\Pi=\sum _{l=1}^2 \Pi^{(lk)}$, and by (\ref{eqn:V-tau-only-matters})
\begin{eqnarray}
\label{eqn:split-in-pi}
\lefteqn{\int _{-\infty}^{+\infty}G_k(k\lambda-\eta)\, d\mu_{\mathcal{T}^{(k)}}(\eta)}\\
&\sim& k\,\sum _{l=1}^2
\int_{1/2}^{\lambda}
\left(V_{\chi_k\,e^{-ikb(\cdot)}}\circ \Pi_k^{(lk)}\right)
\left(x,x+\frac{\mathbf{v}}{\sqrt{k}}\right)\,db.\nonumber
\end{eqnarray}

\begin{lem}
\label{lem:second-term-pi-decays}
There exist $N>0$ such that the following holds.
As
$k\rightarrow +\infty$,
uniformly in $x\in X$, $\|\mathbf{v}\|\le k^{\varpi}$ and $b\in [1/2,+\infty)$ we have
$$\left(V_{\chi_k\,e^{-ikb(\cdot)}}\circ \Pi_k^{(2k)}\right)
\left(x,x+\frac{\mathbf{v}}{\sqrt{k}}\right)=b^N\,O\left(k^{-\infty}\right).
$$
\end{lem}

\textit{Proof.}
By definition,
\begin{eqnarray}
\label{eqn:composition-as-integral}
\lefteqn{\left(V_{\chi_k\,e^{-ikb(\cdot)}}\circ \Pi_k^{(2k)}\right)
\left(x,x+\frac{\mathbf{v}}{\sqrt{k}}\right)}\nonumber
\\
&=&
\int_X\,V_{\chi_k\,e^{-ikb(\cdot)}}\left(x,y\right)\,
\Pi_k^{(2k)}
\left(y,x+\frac{\mathbf{v}}{\sqrt{k}}\right)\,d\mu_X(y);
\end{eqnarray}
it is then enough to prove the claimed estimate for the integrand in (\ref{eqn:composition-as-integral}).

By assumption $Q\ge \mathrm{id}$.
Let $1\le \mu_1\le \mu_2\le \ldots$ be the eigenvalues of $Q$ on $L^2(X)$, repeated according to
multiplicity; thus $\{\lambda_{kj}\}\subseteq \{\mu_\ell\}$.
Let $(\upsilon_\ell)$ be a complete orthonormal system in $L^2(X)$, such that each $\upsilon_\ell$
is an eigenvector of $Q$ with eigenvalue $\mu_\ell$.
Then $\|\upsilon_\ell\|_{\mathcal{C}^0}\le D\,\mu_\ell ^r$ for some fixed $D,r\in \mathbb{R}$,
and there exists
$N>0$ such that $\sum_\ell \mu_\ell^{2r-N}$ converges
(\S 12 of \cite{gs}).

Define
$$
F_k\left(\zeta,x',x''\right)=:\sum _\ell \widehat{\chi}_k
\left(\zeta-\mu_\ell\right)\,\upsilon_\ell\left(x'\right)\,\overline{\upsilon_\ell\left(x''\right)}
\,\,\,\,\,\,\,\,\left(\zeta\in \mathbb{R},\,x',x''\in X\right),
$$
so that
$V_{\chi_k\,e^{-ikb(\cdot)}}\left(x',x''\right)=F_k\left(kb,x',x''\right)$.

By (\ref{eqn:rescaled-rapid-decay}),
for each $N>0$ we have
$\left|\widehat{\chi}(a)\right|<C_N\,k^{(N-1)\xi}\big(1+|a|\big)^{-N}$ for every
$a\in \mathbb{R}$, where $C_N>0$ is constant.
As $\left(1/\mu_\ell\right)+\left|\left(\zeta/\mu_\ell\right)-1\right|\ge 1/\zeta$
for every $\zeta\ge 1$ and $\ell=1,2,\ldots$, we have
\begin{eqnarray*}
\lefteqn{\left|F_k\left(\zeta,x',x''\right)\right|\le C_N\,k^{(N-1)\xi}\,\sum_\ell \big[1+
\left|\zeta-\mu_\ell\right|\big]^{-N}\,\mu_\ell^{2r}}\\
&=&C_N\,k^{(N-1)\xi}\,\sum_\ell \left[\frac{1}{\mu_\ell}+
\left|\frac{\zeta}{\mu_\ell}-1\right|\right]^{-N}\,\mu_\ell^{2r-N}\\
&\le&
C_N\,k^{(N-1)\xi}\,\zeta^N\,\sum_\ell \mu_\ell^{2r-N}\le C'_N\,k^{(N-1)\xi}\,\zeta^N.\nonumber
\end{eqnarray*}
Setting $\zeta=k\,b$ for $b\ge 1/2$, $k\ge 2$ we deduce
\begin{eqnarray}
\label{eqn:polynomial-growth}
\left|V_{\chi\,e^{-ikb(\cdot)}}\left(x',x''\right)\right|
\le C'_N\,k^{(N-1)\xi+N}\,b^N,
\end{eqnarray}
so it suffices to show that the second factor in the integrand in
(\ref{eqn:composition-as-integral}) is $O\left(k^{-\infty}\right)$.

Because $\varrho^{2k}\circ \pi$ is $S^1$-invariant,
$$
\Pi_k^{(2k)}\left(y,x+\frac{\mathbf{v}}{\sqrt{k}}\right)=
\varrho^{(2k)}\big(\pi(y)\big)\,\Pi_k\left(y,x+\frac{\mathbf{v}}{\sqrt{k}}\right).
$$
Let us write $\mathrm{dist}_M=\mathrm{dist}_M\circ (\pi\times \pi):X\times X\rightarrow
\mathbb{R}$.
Then
$\|\mathbf{v}\|\le k^{\varpi}\,\Rightarrow\,\mathrm{dist}_M\left(x,x+\mathbf{v}/\sqrt{k}\right)
\le (5/4)\,k^{\varpi-1/2}$,
hence if $y\in \mathrm{supp}\left(\varrho^{(2k)}\circ \pi\right)$ by (\ref{eqn:bound-on-distance-rho-2})
we have
$
\mathrm{dist}_M\left(y,x+\mathbf{v}/\sqrt{k}\right)\ge
(3/2)\,k^{\varpi-1/2}-(5/4)\,k^{\varpi-1/2}=(1/4)\,k^{\varpi-1/2}$.
By (6.1) of \cite{c},
$\left|\Pi_k\left(y,x+\mathbf{v}/\sqrt{k}\right)\right|\le A\,k^\mathrm{d}\,e^{-B\,k^{\varpi}}$,
for certain positive constants $A,B$.

\hfill Q.E.D.

\begin{rem}
The same argument applies under the hypothesis $\|\mathbf{v}\|\le S\,k^{\varpi}$, for a constant
$S>0$, with the provision that in (\ref{eqn:def-of-varrho-lk}) $\beta$ be replaced with
$\beta/S$.
\end{rem}

\begin{cor}
\label{cor:second-term-pi-decays}
Uniformly for $\|\mathbf{v}\|\lesssim k^{\varpi}$ as
$k\rightarrow +\infty$ we have
$$\int_{1/2}^{\lambda}
\left(V_{\chi_k\,e^{-ikb(\cdot)}}\circ \Pi_k^{(2k)}\right)
\left(x,x+\frac{\mathbf{v}}{\sqrt{k}}\right)\,db=O\left(k^{-\infty}\right).
$$
\end{cor}

We are now reduced to studying the asymptotics of the first summand on the right
hand side of (\ref{eqn:split-in-pi}).
As in the proof of Lemma \ref{lem:second-term-pi-decays},
in so doing by (\ref{eqn:polynomial-growth}) we may modify
$\Pi$ at will as far as the change induced in $\Pi_k\left(y,x+\mathbf{v}/\sqrt{k}\right)$
is rapidly decaying as $k\rightarrow +\infty$.

By \cite{bs}, we can write
$\Pi=\Pi'+\Pi''$, where $\Pi''$ is smoothing and (locally)
the Schwartz kernel of $\Pi'$ is a Fourier integral
\begin{equation}
\label{eqn:pi-as-fio}
\Pi'\left(x',x''\right)=\int_0^{+\infty}e^{it\psi\left(x',x''\right)}\,s\left(t,x',x''\right)\,dt;
\end{equation}
the phase satisfies $\Im\psi\ge 0$ and its Taylor series along the diagonal is determined by the
Hermitian structure, while the amplitude is a classical symbol, $s\left(t,x',x''\right)\sim
\sum_{r=0}^{+\infty}t^{\mathrm{d}-r}\,s_r\left(x',x''\right)$.

Since $\Pi''$ is smoothing,
$\Pi^{''}_k=O\left(k^{-\infty}\right)$, thus
we shall implicitly replace $\Pi$ with $\Pi'$ in (\ref{eqn:V-tau-only-matters}).
In addition, recalling the microlocal structure of $\Pi$, we also have the following
reduction.

\begin{rem}
\label{lem:localize-in-theta}
In Heisenberg local coordinates write $y=x+(\theta,\mathbf{u})$.
On the domain of integration,
$\mathrm{dist}_X\big(x+(\theta,\mathbf{u}),x+\mathbf{v}\big)\ge |\theta|/2$, say.
Therefore,
since $\Pi$ is smoothing away from the diagonal of $X$, we again only lose a rapidly decaying
contribution introducing a cut-off in
$\theta$ which is identically one near $0$ and vanishes for $|\theta|>\varepsilon$
for some small $\varepsilon>0$. We shall also implicitly absorb this cut-off in $s$.
\end{rem}

\begin{rem}
\label{rem:sum-up}
With these reductions understood, in view of Lemmata \ref{lem:irrelevant-interval},
\ref{lem:S''-rapid-decrease} and \ref{lem:second-term-pi-decays},
and of (\ref{eqn:V-tau}), (\ref{eqn:phase-v-tau}), and (\ref{eqn:pi-as-fio}),
we can summarize the previous results on the
asymptotics of $S_{\chi\,e^{-ikb(\cdot)}}^{(k)}$ as follows.

\begin{enumerate}
  \item If $b\not\in (1/2,A_{\varsigma_T}+1)$, then uniformly on $X\times X$
  $$
  \left|S_{\chi_k\,e^{-ikb(\cdot)}}^{(k)}\left(x',x''\right)\right|=O\left(\big(|b|+1\big)^{-\infty}\right)
  \cdot O\left(k^{-\infty}\right).
  $$
  \item If $b\in (1/2,A_{\varsigma_T}+1)$, then uniformly for $\|\mathbf{v}\|\lesssim k^{1/6}$
  \begin{eqnarray}
  \label{eqn:sum-up-asym-S}
  \lefteqn{S_{\chi_k\,e^{-ikb(\cdot)}}^{(k)}\left(x,x+\frac{\mathbf{v}}{\sqrt{k}}\right)=O\left(k^{-\infty}\right)}\\
  &&+\frac{1}{(2\pi)^{2\mathrm{d}+2}}\,
\int_X\int_{-\pi}^\pi\int _{-\epsilon}^\epsilon\int_{\mathbb{R}^{2\mathrm{d}+1}}
\int _0^{+\infty}e^{i\Phi_1}\,A_k\,d\mu_X(y)\,d\vartheta\,d\tau\,d\eta\,dt.\nonumber
  \end{eqnarray}
  \end{enumerate}
  Here $\Phi_1$ and $A$ are as follows. By the $S^1$-invariance of $Q$,
\begin{eqnarray}
\label{eqn:1-st-integral-form-phase}
\Phi_1&=&\varphi\big(\tau,r_\vartheta(x),\eta\big)-y\cdot \eta+t\,\psi\left(y,x+\frac{\mathbf{v}}{\sqrt{k}}\right)
-k\,(b\tau+\vartheta)\\
&=&\big(r_\vartheta(x)-y\big)\cdot \eta+\tau\,q(x,\eta)
+t\,\psi\left(y,x+\frac{\mathbf{v}}{\sqrt{k}}\right)-k\,(b\tau+\vartheta)+
O\left(\tau^2\right)\cdot \|\eta\|;\nonumber
\end{eqnarray}
in Heisenberg local coordinates $r_\vartheta(x)=x+(\vartheta,\mathbf{0})$,
$y=x+(\theta,\mathbf{u})$, and the first summand in (\ref{eqn:1-st-integral-form-phase})
is $-(\theta-\vartheta,\mathbf{u})\cdot \eta$. Also,
\begin{equation}
\label{eqn:1-st-integral-form-amplitude}
A_k=\chi_k(\tau)\,\varrho^{(1k)}(y)\,a\big(\tau,r_\vartheta(x),y,\eta\big)
\,s\left(y,x+\frac{\mathbf{v}}{\sqrt{k}}\right).
\end{equation}
 \end{rem}

\bigskip

We now show that, perhaps after disregarding a rapidly decreasing contribution,
integration in $d\eta$ in (\ref{eqn:sum-up-asym-S})
can be restricted to suitable annuli centered at the
origin, whose radii grow linearly with $k$.
To this end, following \cite{gs}
let $F\in \mathcal{C}^\infty_0(\mathbb{R})$ be identically
equal to $1$ on $(1/C,C)$, for some $C\gg 0$. Write $d\nu$ for the collective integration variables in
(\ref{eqn:sum-up-asym-S}).
If $b\in (1/2,A_{\varsigma_T}+1)$ then
\begin{eqnarray}
\label{eqn:1-st-integral-form-F-compact}
\lefteqn{S_{\chi_k\,e^{-ikb(\cdot)}}^{(k)}\left(x,x+\frac{\mathbf{v}}{\sqrt{k}}\right)}\\
&\sim &\frac{1}{(2\pi)^{2\mathrm{d}+2}}\,
\int_X\int_{-\pi}^\pi\int _{-\epsilon}^\epsilon\int_{\mathbb{R}^{2\mathrm{d}+1}}
\int _0^{+\infty}e^{i\Phi_1}\,F\left(\frac{\|\eta\|}{k}\right)\,A_k\,d\nu\nonumber\\
&&+\frac{1}{(2\pi)^{2\mathrm{d}+2}}\,
\int_X\int_{-\pi}^\pi\int
_{-\epsilon}^\epsilon\int_{\mathbb{R}^{2\mathrm{d}+1}} \int
_0^{+\infty}e^{i\Phi_1}\,\left[1-F\left(\frac{\|\eta\|}{k}\right)\right]\,A_k\,d\nu.\nonumber
\end{eqnarray}

\begin{lem}

Given that $1/2\le b\le A_f+1$,
the latter summand is $O\left(k^{-\infty}\right)$.
\end{lem}
\label{lem:1-st-integral-form-F-compact}
\textit{Proof.} If $F(\|\eta\|/k)\neq 1$, then either
$k\ge C\,\|\eta\|$, or else $\|\eta\|\ge C\,k$. Recall that
$1/2\le b\le A_f+1$.

If $\|\eta\|\ge C\,k$, since $q$ is an elliptic symbol for some $C'>0$ we have
\begin{eqnarray*}
\lefteqn{\big|\partial_\tau \Phi_1\big|=\big|q(x,\eta)-kb\big|}\\
&\ge&C'\,\|\eta\|-kb\ge\frac 12\,C'\,\|\eta\|+\frac 12\,\left(C'C-2 b\right)\,k;\nonumber
\end{eqnarray*}
similarly, if $k\ge C\,\|\eta\|$ then for some $C''>0$
\begin{eqnarray*}
\lefteqn{\big|\partial_\tau \Phi_1\big|=\big|q(x,\eta)-kb\big|}\\
&\ge&k\,b-q(x,\eta)\ge\frac 12\,k -C''\,\|\eta\|\ge\frac 14\,k+\frac 14\,\left(C-4\,C''\right)\,\|\eta\|.\nonumber
\end{eqnarray*}
Therefore, where $F(\|\eta\|/k)\neq 1$ we have
$\big|\partial_\tau \Phi_1\big|\ge c\,k+d\,\|\eta\|$ for some constants $c,\,d>0$.
The claim follows by successively integrating by parts in $d\tau$, since each step introduces a factor
$O\left(k^{\xi-1}\right)$.

\hfill Q.E.D.

 \begin{rem}
\label{rem:choose-C-epsilon}
For future reference, we notice that how large $C$ has to be only depends on
the symbol $q$, and not on the chosen $\epsilon>0$ bounding the size of the
support of the test function $\chi$. In particular, we may assume that
the product $\tau\,\|\eta\|$ is arbitrarily small on the restricted
domain of integration.
\end{rem}

\bigskip

By Lemma
\ref{lem:1-st-integral-form-F-compact}, we need only consider the first
summand in (\ref{eqn:1-st-integral-form-F-compact}).
Recalling (\ref{eqn:split-in-pi}) and
Corollary \ref{cor:second-term-pi-decays},  we get
\begin{eqnarray}
\label{eqn:1-st-integral-form}
\lefteqn{\int _{-\infty}^{+\infty}G_k(k\lambda-\eta)\, d\mu_{\mathcal{T}^{(k)}}(\eta)}\\
&\sim&
\frac{k}{(2\pi)^{2\mathrm{d}+2}}\,
\int_{1/2}^\lambda\,\int_X\int_{-\pi}^\pi\int _{-\epsilon}^\epsilon\int_{\mathbb{R}^{2\mathrm{d}+1}}
\int _0^{+\infty}e^{i\Phi_1}\,F\left(\frac{\|\eta\|}{k}\right)
\,A_k\,db\,d\mu_X(y)\,d\vartheta\,d\tau\,d\eta\,dt.\nonumber
\end{eqnarray}

In Heisenberg local coordinates,
we shall write $y=x+(\theta,\mathbf{u})$ and $d\mu_X(y)=\mathcal{V}(\theta,\mathbf{u})
\,d\theta\,d\mathbf{u}$. With the change of integration variables
$\eta\rightsquigarrow k\,\eta$, $t\rightsquigarrow k\,t$ we get
\begin{eqnarray}
\label{eqn:1-st-integral-form-F-oscillatory}
\lefteqn{\int _{-\infty}^{+\infty}G_k(k\lambda-\eta)\, d\mu_{\mathcal{T}^{(k)}}(\eta)}\\
&\sim&k\,\int_{1/2}^\lambda \left[\left(\frac{k}{2\pi}\right)^{2\mathrm{d}+2}\,
\int_{\mathbb{C}^\mathrm{d}} \int _{-\varepsilon}^{\varepsilon}\int
_{-\pi}^\pi\int_{-\epsilon}^\epsilon\int_{\mathbb{R}^{2\mathrm{d}+1}}
\int_0^{+\infty}e^{ik\Phi_2}\,B_k\, d\widetilde{\nu}\right]\,db,\nonumber
\end{eqnarray}
where
$d\widetilde{\nu}=:\mathcal{V}(\theta,\mathbf{u})\,d\mathbf{u}\,d\theta\,d\vartheta\,
d\tau\,d\eta\,dt$, and
\begin{eqnarray}
\label{eqn:1-st-integral-form-F-oscillatory-phase}
\Phi_2&=:&-(\theta-\vartheta,\mathbf{u})\cdot \eta+\tau\,q(x,\eta)+t\,\psi\left(x+(\theta,\mathbf{u}),
x+\frac{v}{\sqrt{k}}\right)\nonumber\\
&&-(b\tau+\vartheta)+O\left(\tau^2\right)\cdot \|\eta\|,
\end{eqnarray}
\begin{eqnarray}
\label{eqn:1-st-integral-form-F-oscillatory-amplitude}
B_k&=:&F\big(\|\eta\|\big)\,
\chi_k(\tau)\cdot \varrho^{(1k)}\big(x+(\theta,\mathbf{u})\big)\nonumber\\
&&\cdot
a\Big(\tau,r_\vartheta(x),x+(\theta,\mathbf{u}),k\,\eta\Big)\,
s\left(x+(\theta,\mathbf{u}),x+\frac{\mathbf{v}}{\sqrt{k}},k\,t\right).
\end{eqnarray}
Integration in $d\eta$ is over the
ring $1/C\le \|\eta\|\le C$, while in view of the factor $\varrho^{(1k)}\big(x+(\theta,\mathbf{u})\big)$
integration in $d\mathbf{u}$ is over a ball of radius
$O\left(k^{-1/3}\right)$. The inner integral in (\ref{eqn:1-st-integral-form-F-oscillatory})
yields the asymptotics for
$S^{(k)}_{\chi_k \, e^{-ikb(\cdot)}}\left(x,x+\mathbf{v}/\sqrt{k}\right)$.

Let us now set $\eta=r\,\omega$, where $r>0$ and $\omega\in
S^{2\mathrm{d}}$; thus, $\omega=(\omega_0,\mathbf{\omega}_1)\in
\mathbb{R}\times \mathbb{R}^{2\mathrm{d}+1}$,
$\omega_0^2+\|\omega_1\|^2=1$. We have
$d\eta=r^{2\mathrm{d}}\,dr\,d\omega$, and integration in $dr$ is
over $\left(1/C,C\right)$. Thus
(\ref{eqn:1-st-integral-form-F-oscillatory}) may be rewritten:

\begin{eqnarray}
\label{eqn:1-st-integral-form-F-oscillatory-polar}
\lefteqn{\int _{-\infty}^{+\infty}G_k(k\lambda-\eta)\, d\mu_{\mathcal{T}^{(k)}}(\eta)}\\
&\sim&k\,\int_{1/2}^\lambda \left[\left(\frac{k}{2\pi}\right)^{2\mathrm{d}+2}\,
\int_{\mathbb{C}^\mathrm{d}} \int _{-\varepsilon}^\varepsilon\int
_{-\pi}^\pi\int_{-\epsilon}^\epsilon\int_{1/C}^C\int_{S^{2\mathrm{d}}}
\int_0^{+\infty}e^{ik\Phi_3}\,B_k\, d\widehat{\nu}\right]\,db,\nonumber
\end{eqnarray}
where with abuse of language $B_k$ is the amplitude of
(\ref{eqn:1-st-integral-form-F-oscillatory}) with the new variables
inserted,
$d\widehat{\nu}=:r\,^{2\mathrm{d}}\,\mathcal{V}(\theta,\mathbf{u})\,d\mathbf{u}
\,d\theta\,d\vartheta\, d\tau\,dr\,d\omega\,dt$, and
\begin{eqnarray}
\label{eqn:1-st-integral-form-F-oscillatory-phase-polar}
\Phi_3&=:&-r\,(\theta-\vartheta)\,\omega_0-r\,\mathbf{u}\cdot
\omega_1+r\,\tau\,q(x,\omega)+t\,\psi\left(x+(\theta,\mathbf{u}),
x+\frac{\mathbf{v}}{\sqrt{k}}\right)\nonumber\\
&&-(b\tau+\vartheta)+O\left(\tau^2\right)\cdot r.
\end{eqnarray}

The following remark is in order.

\begin{rem}
\label{rem:compact-in-vartheta}
Up to a factor $2\pi$, the
integration $\int _{-\pi}^\pi\,d\vartheta$ is really an integration
over $S^1$. Let $S^1=U_1\cup U_2$ be an open cover, such that
$1\not\in \overline{U}_2$ and $\overline{U}_1\subseteq S^1\setminus
\{-1\}$, and let $\{\rho_1,\rho_2\}$ be a partition of unity
subordinate to it; hence $\rho_1\equiv 1$ in a neighborhood of $1$,
and $\rho_1\equiv 0$ in a neighborhood of $-1$. Then as integration
operators $\int_{S^1}dg= \int_{U_1}\rho_1(g)dg+ \int
_{U_2}\rho_2(g)\,dg$.
With this understanding, we may manipulate
integration in $d\theta$ as if it were compactly supported in
$\theta$, clearly changing the interval of integration in expressing
the second summand; the following computation will show however that
only the first summand contributes non-negligibly.
To keep notation
simple, we shall leave this amendment implicit in the following
discussion.
\end{rem}

\begin{lem}
\label{lem:restrict-on-sphere}
Perhaps after disregarding yet another term
$O\left(k^{-\infty}\right)$, in
(\ref{eqn:1-st-integral-form-F-oscillatory-polar}) integration in
$d\omega$ may be restricted to an appropriate open subset
$S=S_\delta\subseteq S^{2\mathrm{d}}$ of the form $\omega_0>\delta$,
for some $\delta>0$.
\end{lem}

\textit{Proof.}
Since $r\le C$, if $\omega_0<1/(2C)$ we get
$\left|\partial _\vartheta\Phi_3\right|=|r\,\omega_0-1|\ge
1-1/2=1/2$; integration by parts in $d\vartheta$ then proves the
claim with, say, $\delta=1/(3C)$.

\hfill Q.E.D.

\bigskip

Again, this implies introducing a
smooth cut-off identically equal to $1$ in a
neighborhood of the north pole $(1,\mathbf{0})$.
This cut-off will be implicitly absorbed in the amplitude $B_k$;

On $S$,
$\omega_0=\sqrt{1-\|\omega_1\|^2}$, and $\omega_1\in
B_{2\mathrm{d}}\left(\mathbf{0},\sqrt{1-\delta^2}\right)$ is a
system of local coordinates.

\begin{lem}
Perhaps after disregarding a rapidly decreasing contribution,
integration in $dt$ may be restricted to a
compact interval $(1/D,D)$ for some $D\gg 0$.
\end{lem}

\textit{Proof.} The proof will be sketchy, as it
parallels similar arguments in \cite{p}.
As $k\rightarrow +\infty$, uniformly on the given domain of
integration the point
$\left(x+(\theta,\mathbf{u}),x+\mathbf{v}/\sqrt{k}\right)$ is
arbitrarily close to $\big(x+(\theta,\mathbf{0}),x\big)$. Now at the
latter point the differential of $\psi$ is
$\left(e^{i\theta}\,\alpha_{e^{i\theta}x},-e^{-i\theta}\,\alpha_{x}\right)$.
Thus if
$\Upsilon(\theta,\mathbf{u},\mathbf{v})=:\psi\left(x+(\theta,\mathbf{u}),
x+\mathbf{v}/\sqrt{k}\right)$, then for all $k\gg 0$ we have
$3/2\ge \left|\partial_\theta\Upsilon\right|\ge 1/2$.

If $\eta=(\eta_0,\eta_1)\in \mathbb{R}\times
\mathbb{R}^{2\mathrm{d}+1}$ then
$-(\theta-\vartheta,\mathbf{u})\cdot
\eta=(\vartheta-\theta)\,\eta_0+\mathbf{u}\cdot \eta_1$. Since
$|\eta_0|\le \|\eta\|\le C$, we then see from
(\ref{eqn:1-st-integral-form-F-oscillatory-phase}) that
$\left|\partial _\theta\Phi_2\right|\ge (1/2)\,t-C\ge (1/4)\,t+(C/2)$ if
$t>6C$, say. Integrating by parts in $d\theta$ then shows that the
contribution coming from the interval $(6C,+\infty)$ is of rapid
decay.

Similarly, we also deduce from
(\ref{eqn:1-st-integral-form-F-oscillatory-phase}) that
$\left|\partial _\theta\Phi_2\right|\ge r\,\omega_0-(3/2)\,t$.
Since $r\,\omega_0\ge \delta/C$, we get
$\left|\partial _\theta\Phi_2\right|\ge (\delta/C)-(3/2)\,t>\delta/2$
if, say $t<(\delta/3C)$. Therefore, the corresponding contribution
is also $O\left(k^{-\infty}\right)$.

\hfill Q.E.D.

\begin{rem}
\label{rem:compac-in-t}
Thus we can introduce a smooth cut-off in $t$, $\varrho\in \mathcal{C}^\infty_0(\mathbb{R})$,
identically one on $\left(1/D',D'\right)$ and supported in $\left(1/D,D\right)$ for some $0\ll D'\ll D$,
without affecting the asymptotics in (\ref{eqn:1-st-integral-form-F-stainary-phase}).
The cut-off will be henceforth tacitly absorbed into the amplitude, and $\int_0^{+\infty}dt$
replaced with $\int_1^{D}dt$.
\end{rem}

Let us make the change of integration variable
$u\rightsquigarrow u/(r\sqrt{k})$, so that $du\rightsquigarrow
du/\left(r^{2\mathrm{d}}\,k^{\mathrm{d}}\right)$, and integration
in $du$ is now over a ball of radius $O\left(k^{1/6}\right)$.
By (65) of \cite{sz},
\begin{eqnarray*}
\lefteqn{t\,\psi\left(x+\left(\theta,\frac{\mathbf{u}}{r\sqrt{k}}\right),
x+\frac{\mathbf{v}}{\sqrt{k}}\right)}\\
&=&i\,t\,\left[1-e^{i\theta}\right]-\frac{it}{k}\,\psi_2\left(\frac{\mathbf{u}}{r},
\mathbf{v}\right)\,e^{i\theta}
+t\,R_3^\psi\left(\frac{\mathbf{u}}{r\sqrt{k}},\frac{\mathbf{v}}{\sqrt{k}}\right)
\,e^{i\theta}.
\end{eqnarray*}
Inserting this in (\ref{eqn:1-st-integral-form-F-oscillatory-phase-polar}),
(\ref{eqn:1-st-integral-form-F-oscillatory-polar}) may be rewritten
\begin{eqnarray}
\label{eqn:1-st-integral-form-F-stainary-phase}
\lefteqn{\int _{-\infty}^{+\infty}G_k(k\lambda-\eta)\, d\mu_{\mathcal{T}^{(k)}}(\eta)\sim k}\\
&&\cdot
\int_{1/2}^\lambda \left\{\frac{k^{\mathrm{d}+2}}{(2\pi)^{2\mathrm{d}+2}}\int_{\mathbb{C}^\mathrm{d}}\int_{S^{2\mathrm{d}}}e^{-i\sqrt{k}\,\mathbf{u}\cdot \omega_1}
\left[\int_0^{+\infty}\int _{-\varepsilon}^\varepsilon\int_{-\pi}^\pi
\int_{1/C}^C\int_{-\epsilon}^\epsilon
e^{ik\Phi}\,J_k\, d\gamma\right]\,d\mathbf{u}\,d\omega\right\}\,db,\nonumber
\end{eqnarray}
where $\mathbf{u}\cdot \omega_1$ is the standard Euclidean pairing
in $\mathbb{R}^{2\mathrm{d}}$,
$d\gamma=:dt\,d\theta\,d\vartheta\,dr\,d\tau$,
\begin{eqnarray}
\label{eqn:1-st-integral-form-stationary-phase}
\Phi&=:&i\,t\,\left[1-e^{i\theta}\right]-r\,(\theta-\vartheta)\,\omega_0+r\,\tau\,q(x,\omega)
-(b\tau+\vartheta)\nonumber\\
&&+O\left(\tau^2\right)\cdot r,
\end{eqnarray}
and
\begin{equation}
\label{eqn:1-st-integral-form-stationary-amplitude}
J_k=:e^{t\psi_2\left(\mathbf{u}/r,\mathbf{v}\right)\,e^{i\theta}+itkR_3^\psi\left(\mathbf{u}
/(r\sqrt{k}),\mathbf{v}/\sqrt{k}\right)\,e^{i\theta}}\,B_k.
\end{equation}
We have $\Im\Phi\ge 0$ and
$|J_k|\le C'\,k^{\mathrm{d}}\,t^\mathrm{d}\,e^{-a\|\mathbf{v}-\mathbf{u}\|^2}$ on the domain of
integration, for an appropriate $a>0$; also,
$tkR_3^\psi\left(\mathbf{u}
/(r\sqrt{k}),\mathbf{v}/\sqrt{k}\right)\,e^{i\theta}$ remains bounded
for $\|\mathbf{u}\|,\|\mathbf{v}\|=O\left(k^{1/6}\right)$.
The expression within $\{\,\}$ in (\ref{eqn:1-st-integral-form-F-stainary-phase})
yields the asymptotics for $S^{(k)}_{\chi_k\,e^{-ikb(\cdot)}}\left(x,x+\mathbf{v}/\sqrt{k}\right)$.

\begin{rem}
\label{rem:leading-order-term}
By construction, $J_k=\chi\left(k^\xi\right)\,L_k$, where $L_k$ is the product of
the rescaled amplitude
$a\left(\tau,r_\vartheta(x),x+\big(\theta,\mathbf{u}/(r\sqrt{k}\big),k\eta\right)$ of
$V$, the rescaled amplitude $s\left(x+\big(\theta,\mathbf{u}/(r\sqrt{k}\big),x+\mathbf{v}/\sqrt{k},kt\right)$
of $\Pi$,
the various cut-offs introduced, the exponential factor
in (\ref{eqn:1-st-integral-form-stationary-amplitude}), and the local coordinate expression
for the Riemannian density of $X$. Now $s$ and $a$ are semiclassical
symbols, and as such they admit asymptotic expansions in descending powers of $k$.
On the other hand, all the factors involved may be Taylor expanded in the
arguments $\mathbf{u}/(r\sqrt{k})$ and $\mathbf{v}/\sqrt{k}$, so that $L_k$
has an asymptotic expansion in descending powers of $k^{-1/2}$. To determine the leading
order term of the asymptotic expansion of $L_k$ where the cut-offs are identically $1$,
we remark that the leading order terms coming from $a$ and $s$
are, respectively, $1/\mathcal{V}\big(x+\big(\theta,\mathbf{0})\big)$ and $(k/\pi)^\mathrm{d} \,t^{\mathrm{d}}$;
the first factor cancels with the $\mathcal{V}$ term in the local coordinate expression for the
density. As a result, the leading order term
is $(k/\pi)^\mathrm{d} \,t^{\mathrm{d}}\,e^{t\,e^{i\theta}\psi_2\left(\mathbf{u}/r,\mathbf{v}\right)}$.
\end{rem}

Now we notice that by construction $\Phi$ does not depend on $\mathbf{u}$; therefore, the
exponential $e^{ik\Phi}$ behaves
like a constant under differentiation by the $\mathbf{u}$ variables. Let us view the right
hand side of (\ref{eqn:1-st-integral-form-F-stainary-phase}) in its dependence on $\mathbf{u}$
as an oscillatory integral
with phase $\Psi=-i\sqrt{k}\,\mathbf{u}\cdot \omega$. Fix a sufficiently small $t_0>0$ .
Since $\nabla_{\mathbf{u}}\Psi=-i\,\omega_1$,
on the locus $T_k\subseteq S_\delta$ where $\|\omega_1\|>(t_0/2)\,k^{-\xi}$
(equivalently, $\omega_0<\sqrt{1-(t_0^2/4)\,k^{-2\xi}}$)
we can successively integrate by parts in $\mathbf{u}$, and at each step we introduce a factor
$k^{\xi-1/2}$. Therefore, for any $t_0>0$
the contribution to (\ref{eqn:1-st-integral-form-F-stainary-phase})
coming from $T_k$
is $O\left(k^{-\infty}\right)$.
Perhaps after disregarding a rapidly decreasing contribution,
we may thus introduce an appropriate partition of unity on $S_\delta$, of the form
$\gamma_{kj}\left(k^\xi\,\omega_1\right)$, and restrict
integration in $\omega_1$ to a progressively shrinking
open neighborhood of $(1,\mathbf{0})$ in $S_\delta$, of the form
$U_k=\left\{\omega_1:\|\omega_1\|<t_0\,k^{-\xi}\right\}$.
We shall now proceed leaving this further cut-off implicit.

This reduction has the following consequence.
Since $q(x,\omega)/\omega_0$ equals $\varsigma_T(m)$ at  $(1,\mathbf{0})$,
and integration only involves points $\omega\in S^{2\mathrm{d}}$
in an $t_0\,k^{-\xi}$-neighborhood of $(1,\mathbf{0})$, on $T_k$ we have
$\big |q(x,\omega)/\omega_0-\varsigma_T(m)\big|<D\,t_0\,k^{-\xi}$
for some fixed $D>0$ and all $k\gg 0$.

There is one further similar reduction that we can make in the variable $r$.
Set
$R_k=:\big\{r:\big|r-1/\omega_0\big|>(t_0/2)k^{-\xi}\big\}$.
Since $\partial_\vartheta \Phi=r\omega_0-1$, there exists $\beta>0$ such that
$\big|\partial_\vartheta \Phi\big|\ge \beta\,k^{-\xi}$ on $R_k$. Successively integrating by parts
in $\vartheta$, we introduce at each step a factor $k^{\xi-1}$; therefore, the contribution
of $R_k$ to the asymptotics of (\ref{eqn:1-st-integral-form-F-stainary-phase})
is $O\left(k^{-\infty}\right)$.
Summing up, upon introducing a cut-off of the form $\sigma\left(k^{\xi}\big(r-1/\omega_0\big)\right)$,
we may reduce integration over $r$ to the interval $S_k$ where
$\big|r-1/\omega_0\big|<t_0k^{-\xi}$.

We are now ready to prove:

\begin{prop}
\label{prop:rapid-decay-case}
Uniformly for
$\big|b-\varsigma_T(m)\big|>e_k$ and
$\|\mathbf{v}\|\lesssim k^{\varpi}$, as $k\rightarrow +\infty$
we have
$$
S_{\chi_k\,e^{-ikb(\cdot)}}^{(k)}\left(x,x+\frac{\mathbf{v}}{\sqrt{k}}\right)=O\left(k^{-\infty}\right).
$$
\end{prop}

\textit{Proof.} Let us consider the inner integral in
(\ref{eqn:1-st-integral-form-F-stainary-phase}) as an oscillatory integral in $\tau$.
We have $\partial_\tau\Phi=r\,q(x,\omega)-b+O(\tau)\cdot r$, and
$|\tau|\le \epsilon\,k^{-\xi}$ on the support of $\chi_k$; therefore,
for some constant $C'>0$,
\begin{eqnarray*}
\lefteqn{\big|\partial_\tau\Phi\big|}\\
&=&\left|\left(r-\frac{1}{\omega_0}\right)\,q(x,\omega)
+\left(\frac{q(x,\omega)}{\omega_0}-\varsigma_T(m)\right)+\big(\varsigma_T(m)-b\big)+O(\tau)\cdot r\right|\\
&\ge&\big|\varsigma_T(m)-b\big|-\left|\frac{q(x,\omega)}{\omega_0}-\varsigma_T(m)\right|-
\left|\left(r-\frac{1}{\omega_0}\right)\,q(x,\omega)\right|-C'\,\epsilon\,k^{-\xi}\\
&\ge&c\,k^{-\xi}-C'\,t_0\,k^{-\xi}-C'\,\epsilon\,k^{-\xi}\ge \frac 12\,c\,k^{-\xi},
\end{eqnarray*}
if $\epsilon$ and $t_0$ have been chosen sufficiently small.

In view of the factor $\chi_k(\tau)=\chi\left(k^\xi \,\tau\right)$,
we conclude that successive partial integrations by parts in $\tau$ introduce at each step a factor
$k^{2\xi-1}$; the statement follows since by hypothesis $\xi<1/2$.

\hfill Q.E.D.

\begin{cor}
\label{cor:rapid-decay-case}
If $\lambda<\varsigma_T(m)$, then
$\int _{-\infty}^{+\infty}G_k(k\lambda-\eta)\, d\mu_{\mathcal{T}^{(k)}}(\eta)=
O\left(k^{-\infty}\right)$ as $k\rightarrow +\infty$.
\end{cor}

Before giving the second computation of
$\int _{-\infty}^{+\infty}G(k\lambda-\eta)\, d\mu_{\mathcal{T}^{(k)}}(\eta)$,
we need to establish some preliminary results.
By Remark \ref{rem:sum-up}
and Proposiion \ref{prop:rapid-decay-case} we have:

\begin{lem}
\label{lem:direct-estimate-on-S}
Fix $C>0$, and define
$\gamma_k(C)=\sup\left\{\gamma_k^+(C),\gamma_k^-(C)\right\}$, where
$$
\gamma_k^+(C)=:\sup\left\{\left|S^{(k)}_{\chi e^{-ik\eta}}(x,x)\right|:\eta\in \left[\varsigma_T(m)+C\,k^{-\xi},+\infty\right)\right\},
$$
and
$$
\gamma_k^-(C)=:\max\left\{\left|S^{(k)}_{\chi e^{-ik\eta}}(x,x)\right|:\eta\in \left(-\infty,\varsigma_T(m)-C\,k^{-\xi}\right]\right\}.
$$
Then $\gamma _k(C)=O\left(k^{-\infty}\right)$ as $k\rightarrow +\infty$.
\end{lem}

At this stage, it is convenient to
make a more specific choice of $\chi$. Namely, choose
first $\psi\in \mathcal{C}^\infty_0(-\epsilon/2,\epsilon/2)$ such that
$\psi\ge 0$, $\psi(t)>0$ for all $t\in (-\epsilon/3,\epsilon/3)$,
$\psi(t)=\psi(-t)$ for all $t\in \mathbb{R}$, $\|\psi\|_{L^2}=1$.
Then the Fourier transform $\widehat{\psi}$ is real. Let
$\chi=:\psi*\psi$; then $\chi\in \mathcal{C}^\infty_0(-\epsilon,\epsilon)$,
$\chi\ge 0$, $\chi(t)>0$ if $t\in (-2\epsilon/3,2\epsilon/3)$, $\chi(0)=\|\psi\|_{L^2}^2=1$.
Furthermore, $\widehat{\chi}=\widehat{\psi}\cdot \widehat{\psi}\ge 0$, and
\begin{equation}
\label{eqn:chi-hat-0}
\widehat{\chi}(0)=\widehat{\psi}(0)^2=\left(\int _{-\infty}^{+\infty}
\psi(\tau)\,d\tau\right)^2=\|\psi\|_{L^1}^2>0.
\end{equation}

\noindent
By (\ref{eqn:chi-hat-0}), there exists $\delta>0$ such that
\begin{equation}
\label{eqn:role-of-delta}
|\lambda|<\delta\,\,\,\,\Rightarrow\,\,\,\,\widehat{\chi}(\lambda)\ge \|\psi\|_{L^1}^2/2.
\end{equation}
Passing to the rescaled function $\chi_k$, (\ref{eqn:role-of-delta}) implies
\begin{equation}
\label{eqn:role-of-delta-k}
|\lambda|<\delta\,k^\xi\,\,\,\,\Rightarrow\,\,\,\,\widehat{\chi}_k(\lambda)\ge \frac{1}{2k^\xi}\,\|\psi\|_{L^1}^2.
\end{equation}

Consequently, if $c\in \mathbb{R}$ and $0\le \delta'\le \delta$, then
\begin{eqnarray}
\label{eqn:bound-by-S}
\lefteqn{\mathcal{T}_k\left(kc+k^\xi\delta',x,x\right)-\mathcal{T}_k(kc,x,x)}\\
&=&\sum_{kc<\lambda_{kj}\le kc+k^\xi\delta'}\left|e_{kj}(x)\right|^2
\le \frac{2k^\xi}{\|\psi\|^2_{L^1}}\sum_{kc<\lambda_{kj}\le kc+\delta'}
\widehat{\chi}_k\big(kc-\lambda_{kj}\big)\,\left|e_{kj}(x)\right|^2\nonumber\\
&\le&\frac{2k^\xi}{\|\psi\|^2_{L^1}}\sum_{j=1}^{N_k}
\widehat{\chi}_k\big(kc-\lambda_{kj}\big)\,\left|e_{kj}(x)\right|^2=
\frac{2k^\xi}{\|\psi\|^2_{L^1}}\,S^{(k)}_{\chi_k e^{-ikc(\cdot)}}(x,x).\nonumber
\end{eqnarray}
A similar estimate holds if $-\delta\le \delta'\le 0$.

\begin{lem}
\label{lem:estimate-techincal-T}
Fix $\ell>0$. Then, uniformly for
$\big|\lambda-\varsigma_T(m)\big|\ge \ell\,k^{-\xi}$, we have
$$
\int_{-\infty}^{+\infty}\Big[\mathcal{T}_k\big(k\lambda+k^{\xi}\,b,x,x\big)-
\mathcal{T}_k\big(k\lambda,x,x\big)\Big]\,\widehat{\chi}(b)\,db=O\left(k^{-\infty}\right)
$$
as $k\rightarrow +\infty$.
\end{lem}

\textit{Proof.}
Let us prove the statement for $\lambda>\varsigma_T(m)$; the case
$\lambda<\varsigma_T(m)$
is similar.
In the present proof,
$\mathcal{T}_k(r)=:\mathcal{T}_k(r,x,x)$ for $r\in \mathbb{R}$,
$k\in \mathbb{N}$.

For $k=1,2,\ldots$ let
$J_1^{(k)}=:\big(-\infty,-(\ell/3)\,k^{1-2\xi}\big)$, $J_2^{(k)}=:(-(\ell/2)\,k^{1-2\xi},+\infty)$.
Thus $\mathcal{J}^{(k)}=:\left\{J_1^{(k)},J_2^{(k)}\right\}$ is an open cover
of $\mathbb{R}$. Let $\{\gamma_j\}$ be a partition of unity on $\mathbb{R}$
subordinate to $\mathcal{J}^{(1)}$. Define
$\gamma_j^{(k)}(b)=:\gamma_j\left(b/k^{1-2\xi}\right)$ ($b\in \mathbb{R}$); then
$\{\gamma_j^{(k)}\}$ is a partition of unity subordinate to $\mathcal{J}^{(k)}$.
For every $k=1,2,\ldots$,
\begin{eqnarray}
\label{eqn:split-total-int-T}
\lefteqn{\int_{-\infty}^{+\infty}\left[\mathcal{T}_k\left(k\lambda+k^\xi b\right)-
\mathcal{T}_k (k\lambda)\right]\,\widehat{\chi}(b)\,db}\nonumber\\
&=&\sum_{j=1}^2\,\int_{-\infty}^{+\infty}\left[\mathcal{T}_k\left(k\lambda+k^\xi b\right)-
\mathcal{T}_k (k\lambda)\right]\,\gamma_j^{(k)}(b)\,\widehat{\chi}(b)\,db.
\end{eqnarray}

Let us first estimate the summand with $j=1$.
On the support of $\gamma_1^{(k)}$ we have $|b|>(\ell/3)k^{1-2\xi}$, hence if $k\gg 0$ then
$|b|\ge |b|/2+(\ell/6)k^{1-2\xi}$.
Since $\widehat{\chi}$ is of rapid decrease, for every integer $N>0$
there exists a constant $C_{2N}>0$
such that
\begin{eqnarray*}
\left|\widehat{\chi}(b)\right|&<&C_{2N}\,\left[\left(\frac 12\,|b|+1\right)+\frac 16\,\ell\,k^{1-2\xi}\right]^{-2N}\\
&\le&\left(\frac 6\ell\right)^N\,C_{2N}\,{2N \choose N}^{-1}\,\left(\frac 12\,|b|+1\right)^{-N}
\,k^{-N(1-2\xi)}.
\end{eqnarray*}
Hence for every $N>0$ there exists a constant $C_N'>0$ such that on the support of
$\gamma_1^{(k)}$ the integrand is bounded by
$C_N'\,\big(|b|/2+1\big)^{-N}\,k^{\mathrm{d}-N}$, and so the first summand
is $O\left(k^{-\infty}\right)$.

Let us now consider the summand with $j=2$ in (\ref{eqn:split-total-int-T}). Let $\delta$
be as in (\ref{eqn:role-of-delta}). Suppose first $b\ge 0$, and write
$b=\delta\,\lfloor b/\delta\rfloor+\delta'$, where $0\le \delta'<\delta$.
Then
\begin{eqnarray}
\label{eqn:split-in-delta}
\lefteqn{\mathcal{T}_k\left(k\lambda+k^\xi b\right)-\mathcal{T}_k(k\lambda)}\\
&=&\left[\mathcal{T}_k\left(k\lambda+k^\xi b\right)-\mathcal{T}_k\left(k\lambda+k^\xi \delta\,\left\lfloor \frac{b}{\delta}\right\rfloor\right)\right]
+\left[\mathcal{T}_k\left(k\lambda+k^\xi \delta\,\left\lfloor \frac{b}{\delta}\right\rfloor\right)-\mathcal{T}_k(k\lambda)\right].
\nonumber
\end{eqnarray}
We apply (\ref{eqn:bound-by-S}) with $c=c_k=:\lambda+\left(\delta/k^{1-\xi}\right)\,\lfloor b/\delta\rfloor$:
\begin{eqnarray}
\label{eqn:bound-by-S-1st-term}
\lefteqn{\mathcal{T}_k\left(k\lambda+k^\xi b\right)
-\mathcal{T}_k\left(k\lambda+k^\xi\delta\,\left\lfloor \frac{b}{\delta}\right\rfloor\right)}\\
&=&\mathcal{T}_k\left(k\,c_k+
k^\xi\delta'\right)-
\mathcal{T}_k\left(k\,c_k\right)\le\frac{2k^\xi}{\|\psi\|^2_{L^1}}\,S^{(k)}_{\chi e^{-ik\,c_k
(\cdot)}}(x,x)\le \frac{2k^\xi}{\|\psi\|^2_{L^1}}\,\gamma_k(\ell),\nonumber
\end{eqnarray}
where $\gamma_k$ is as in Lemma \ref{lem:direct-estimate-on-S}.

Next, setting $b_{kj}=:\lambda+\left(j/k^{1-\xi}\right)\,\delta$,
\begin{eqnarray}
\label{eqn:bound-by-S-2nd-term}
\lefteqn{\mathcal{T}_k\left(k\lambda+k^\xi\left\lfloor \frac{b}{\delta}\right\rfloor\,\delta\right)-\mathcal{T}_k(k\lambda)}\\
&=&\sum_{j=0}^{\lfloor b/\delta\rfloor-1}
\Big[\mathcal{T}_k\left(k\lambda+(j+1)\,k^\xi\delta\right)
-\mathcal{T}_k\left(k\lambda+j\,k^\xi\delta\right)\Big]\nonumber\\
&=&\sum_{j=0}^{\lfloor b/\delta\rfloor-1}
\Big[\mathcal{T}_k\left(k\,b_{kj}+k^\xi\delta\right)
-\mathcal{T}_k\left(k\,b_{kj}\right)\Big]\nonumber\\
&\le&\frac{2k^\xi}{\|\psi\|^2_{L^1}}\,\sum_{j=0}^{\lfloor b/\delta\rfloor-1}S^{(k)}_{\chi e^{-ik\,b_{kj}
(\cdot)}}(x,x)\le \frac{2k^\xi}{\|\psi\|^2_{L^1}}\,\left\lfloor \frac{b}{\delta}\right\rfloor\,\gamma_k(\ell).
\end{eqnarray}

Similarly, if $-k^{1-2\xi}(\ell/2)\le b\le 0$ then writing $b=-b'$, and $b'=\lfloor b'/\delta\rfloor\,\delta +\delta'$
with $0\le \delta'<\delta$ the analogue of (\ref{eqn:split-in-delta}) is
\begin{eqnarray}
\label{eqn:split-in-delta-bneg}
\left|\mathcal{T}_k\left(k\lambda+k^\xi b\right)-
\mathcal{T}_k\big(k\lambda\big)\right|&=&\mathcal{T}_k\big(k\lambda\big)-\mathcal{T}_k\left(k\lambda-k^\xi b'\right)\\
&=&\left[\mathcal{T}_k\big(k\lambda\big)-\mathcal{T}_k\left(k\lambda-
k^\xi\left\lfloor \frac{b'}{\delta}\right\rfloor\,\delta\right)\right]\nonumber\\
&&+
\left[\mathcal{T}_k\left(k\lambda-k^\xi
\left\lfloor \frac{b'}{\delta}\right\rfloor\,\delta\right)-\mathcal{T}_k\left(k\lambda-
k^\xi \left\lfloor \frac{b'}{\delta}\right\rfloor\,\delta-k^\xi\delta'\right)\right].\nonumber
\end{eqnarray}

Let us set $r_k=:\lambda-k^{\xi-1}\left\lfloor b'/\delta\right\rfloor \,\delta$;
since $0\le b'\le (\ell/2)\,k^{1-2\xi}$ we have
\begin{equation*}
\left|r_k
-\varsigma_T(m)\right|\ge
\big|\lambda-\varsigma_T(m)\big|-\frac{1}{k^{1-\xi}}\left\lfloor \frac{b'}{\delta}\right\rfloor \,\delta
\ge \ell\,k^{-\xi}-\frac 12\,\ell\,k^{-\xi}=\frac 12\,\ell\,k^{-\xi}.
\end{equation*}

Therefore, in view of (\ref{eqn:bound-by-S}) we have
\begin{eqnarray*}
\lefteqn{\mathcal{T}_k\left(k\lambda-k^\xi
\left\lfloor \frac{b'}{\delta}\right\rfloor\,\delta\right)-\mathcal{T}_k\left(k\lambda-
k^\xi \left\lfloor \frac{b'}{\delta}\right\rfloor\,\delta-k^\xi\delta'\right)}\\
&=&\mathcal{T}_k\left(k\,r_k\right)-\mathcal{T}_k\left(k\,r_k-k^\xi\delta'\right)\le
\frac{2k^\xi}{\|\psi\|^2_{L^1}}\,S^{(k)}_{\chi e^{-ik\,r_k
(\cdot)}}(x,x)\le \frac{2k^\xi}{\|\psi\|^2_{L^1}}\,\gamma_k\big(\ell/2\big).
\end{eqnarray*}

If furthermore $r_{kj}=:\lambda-\left(j/k^{1-\xi}\right)\,\delta$, then
$$
r_{kj}\ge
\lambda-\left( b'/k^{1-\xi}\right)\ge \lambda-(\ell/2)\,k^{-\xi}\ge \varsigma_T(m)+(\ell/2)\,k^{-\xi}
$$
for $0\le b'\le (\ell/2)\,k^{1-2\xi}$ and $0\le j\le \lfloor b'/\delta\rfloor$.
Therefore, as in (\ref{eqn:bound-by-S-2nd-term}),
\begin{eqnarray*}
\mathcal{T}_k\big(k\lambda\big)-\mathcal{T}_k\left(k\lambda-
k^\xi\left\lfloor \frac{b'}{\delta}\right\rfloor\,\delta\right)&=&\sum_{j=0}^{\lfloor b'/\delta\rfloor-1}\Big[\mathcal{T}_k\left(k\lambda-k^\xi j\,\delta\right)-
\mathcal{T}_k\left(k\lambda-k^\xi (j+1)\,\delta\right)\Big]\\
&=&\sum_{j=0}^{\lfloor b'/\delta\rfloor-1}\Big[\mathcal{T}_k\left(k\,r_{kj}\right)-
\mathcal{T}_k\left(k\,r_{kj}-k^\xi\delta\right)\Big]\\
&\le&\frac{2k^\xi}{\|\psi\|^2_{L^1}}\,\sum_{j=0}^{\lfloor b'/\delta\rfloor-1}
S^{(k)}_{\chi e^{-ik\,r_{kj}
(\cdot)}}(x,x)
\le C\,k^\xi\,|b|\,\gamma_k\big(\ell /2\big).
\end{eqnarray*}
%
%

Summing up, on $J_2^{(k)}$
\begin{equation*}
\left|\Big[\mathcal{T}_k\big(k\lambda+k^\xi b\big)-
\mathcal{T}_k\big(k\lambda\big)\Big]\,\widehat{\chi}(b)\right|\le
C\,k^\xi(1+|b|)\,\widehat{\chi}(b)\,\gamma_k(\ell/2)
\end{equation*}
for some constant $C>0$;
we conclude from Lemma \ref{lem:direct-estimate-on-S} that the summand with $j=2$ in
(\ref{eqn:split-total-int-T}) is also $O\left(k^{-\infty}\right)$.

\hfill Q.E.D.

\begin{cor}
\label{cor:estimate-techincal-T}
Fix $\ell>0$. Then, uniformly in $\big(\lambda,x_1,x_2\big)\in \mathbb{R}\times X\times X$
such that $\big|\lambda-\varsigma_T\left(x_j\right)\big|\ge \ell \,k^{-\xi}$
for $j=1,2$,
we have
$$
\int_{-\infty}^{+\infty}\Big[\mathcal{T}_k\left(k\lambda+k^\xi b,x_1,x_2\right)-
\mathcal{T}_k\left(k\lambda,x_1,x_2\right)\Big]\,\widehat{\chi}(b)\,db=O\left(k^{-\infty}\right).
$$
\end{cor}

\textit{Proof.} This is really a consequence of the proof of Lemma \ref{lem:estimate-techincal-T}.
Suppose $b\ge 0$ (a similar estimate holds for $b\le 0$).
Then
\begin{eqnarray*}
\lefteqn{\Big|\mathcal{T}_k\left(k\lambda+k^\xi b,x_1,x_2\right)-
\mathcal{T}_k\left(k\lambda,x_1,x_2\right)\Big|=\left|\sum_{j:k\lambda<\lambda_{kj}\le k\lambda+k^\xi
b}e_{kj}\left(x_1\right)\,
\overline{e_{kj}\left(x_2\right)}\right|}\\
&\le& \sqrt{\sum_{j:k\lambda<\lambda_{kj}\le k\lambda+k^\xi b}\left|e_{kj}\left(x_1\right)\right|^2}
\cdot \sqrt{\sum_{j:k\lambda<\lambda_{kj}\le k\lambda+k^\xi b}\left|e_{kj}\left(x_2\right)\right|^2}\\
&=&\sqrt{\Big[\mathcal{T}_k\left(k\lambda+k^\xi b,x_1,x_1\right)-
\mathcal{T}_k\left(k\lambda,x_1,x_1\right)\Big]}\cdot \sqrt{\Big[\mathcal{T}_k\left(k\lambda+k^\xi b,x_2,x_2\right)-
\mathcal{T}_k\left(k\lambda,x_2,x_2\right)\Big]}.
\end{eqnarray*}
The statement follows from this and the bounds established in the proof of Lemma \ref{lem:estimate-techincal-T}.

\hfill Q.E.D.

\bigskip

Now we can give the second estimate of
$\int _{-\infty}^{+\infty}G_k(k\lambda-\eta)\, d\mu_{\mathcal{T}^{(k)}}(\eta)$.
Let $H$ denote the Heaviside function.
Recall that $\widehat{\chi}_k(b)=(1/k^\xi)\,\widehat{\chi}\left(b/k^\xi\right)$.
Then
\begin{eqnarray}
\label{eqn:2-nd-computation-key-integral}
\lefteqn{\int _{-\infty}^{+\infty}G_k(k\lambda-\eta)\, d\mu_{\mathcal{T}^{(k)}}(\eta)}\nonumber\\
&=&\sum_{j=1}^{N_k}G_k(k\lambda-\lambda_{kj})\,
e_{kj}(x)\cdot \overline{e_{kj}\left(x+\frac{\mathbf{v}}{\sqrt{k}}\right)}
\nonumber\\
&=&\sum_{j=1}^{N_k}\left(\int_{-\infty}^{k\lambda-\lambda_{kj}}\widehat{\chi}_k(b)\,db\right)\,
e_{kj}(x)\cdot \overline{e_{kj}\left(x+\frac{\mathbf{v}}{\sqrt{k}}\right)}\nonumber\\
&=&\sum_{j=1}^{N_k}\left(\int_{-\infty}^{+\infty}H\big(k\lambda-\lambda_{kj}-b\big)\,
\widehat{\chi}_k(b)\,db\right)\,
e_{kj}(x)\cdot \overline{e_{kj}\left(x+\frac{\mathbf{v}}{\sqrt{k}}\right)}\nonumber\\
&=&\int_{-\infty}^{+\infty}\left[\sum_{j=1}^{N_k}H\big(k\lambda-\lambda_{kj}-b\big)\,
e_{kj}(x)\cdot \overline{e_{kj}\left(x+\frac{\mathbf{v}}{\sqrt{k}}\right)}\right]\,
\widehat{\chi}_k(b)\,db\nonumber\\
&=&\int_{-\infty}^{+\infty}\mathcal{T}_k\left(
k\lambda-b,
x,x+\frac{\mathbf{v}}{\sqrt{k}}\right)\,
\widehat{\chi}_k(b)\,db\nonumber\\
&=&\mathcal{T}_k\left(
k\lambda,
x,x+\frac{\mathbf{v}}{\sqrt{k}}\right)\,\int_{-\infty}^{+\infty}\widehat{\chi}_k(b)\,db\nonumber\\
&&+\int_{-\infty}^{+\infty}\left[\mathcal{T}_k\left(
k\lambda-b,
x,x+\frac{\mathbf{v}}{\sqrt{k}}\right)-\mathcal{T}_k\left(
k\lambda,
x,x+\frac{\mathbf{v}}{\sqrt{k}}\right)\right]\,\widehat{\chi}_k(b)\,db\nonumber\\
&=&2\pi\,\mathcal{T}_k\left(
k\lambda,
x,x+\frac{\mathbf{v}}{\sqrt{k}}\right)\nonumber\\
&&+\int_{-\infty}^{+\infty}\left[\mathcal{T}_k\left(
k\lambda-k^\xi\,b,
x,x+\frac{\mathbf{v}}{\sqrt{k}}\right)-\mathcal{T}_k\left(
k\lambda,
x,x+\frac{\mathbf{v}}{\sqrt{k}}\right)\right]\,\widehat{\chi}(b)\,db\nonumber\\
&=&2\pi\,\mathcal{T}_k\left(
k\lambda,
x,x+\frac{\mathbf{v}}{\sqrt{k}}\right)+O\left(k^{-\infty}\right),
\end{eqnarray}
in view of Corollary \ref{cor:estimate-techincal-T}.

To complete the proof of Theorem \ref{thm:main} we need only compare Corollary \ref{cor:rapid-decay-case}
and (\ref{eqn:2-nd-computation-key-integral}).

\hfill Q.E.D.

\section{Proof of Corollary \ref{cor:asympt-dim}.}

Since it is $S^1$-invariant,
the diagonal restriction $\mathcal{T}_k\left(\eta,x,x\right)$
may regarded as defined on $M$; let us set
$\mathfrak{T}_k\left(\eta,m\right)=:\mathcal{T}_k\left(\eta,x,x\right)$
if $m\in M$ and $\pi\left(x\right)=m$.
Clearly,
\begin{equation*}
\dim\left(V_{k\lambda}^{(k)}\right)=\int_M \mathfrak{T}_k\left(k\lambda,m\right)
\,dV_M(m).
\end{equation*}

For any sufficiently small $\epsilon>0$, consider the disjoint union
$M=M_1^{(\epsilon)}\cup M_2^{(\epsilon)}\cup M_3^{(\epsilon)}$, where
$M_1^{(\epsilon)}=:M_{<\lambda-\epsilon}$, $M_2^{(\epsilon)}=:
M_{<\lambda+\epsilon}\cap M_{<\lambda-\epsilon}^c$, and $M_3^{(\epsilon)}=:
M_{<\lambda+\epsilon}^c$.
By Theorem \ref{thm:main}
we have $\mathfrak{T}_k\left(k\lambda,m\right)=(k/\pi)^{\mathrm{d}}+
O\left(k^{\mathrm{d}-1}\right)$ uniformly on $M_1^{(\epsilon)}$, and
$\mathfrak{T}_k\left(k\lambda,m\right)=
O\left(k^{-\infty}\right)$ uniformly on $M_3^{(\epsilon)}$.
Since $\lambda$ is a regular value of $\varsigma_T$,
$M_2^{(\epsilon)}$ is contained in a $(a\,\epsilon)$-neighborhood of
$M_\lambda$, for some fixed $a>0$; therefore, its volume is
$O(\epsilon)$. By the same token,
$\mathrm{vol}\big(M_{<\lambda}\big)-\mathrm{vol}\big(M_1^{(\epsilon)}\big)=O(\epsilon)$.

Thus
\begin{eqnarray*}
\int_M\,\mathfrak{T}_k\left(k\lambda,m\right)
\,dV_M(m)
&=&\sum _{j=1}^3\int_{M_j^{(\epsilon)}}\,\mathfrak{T}_k\left(k\lambda,m\right)
\,dV_M(m)\\
&=&\mathrm{vol}\big(M_{<\lambda}\big)\,\left(\frac{k}{\pi}\right)^{\mathrm{d}}+O(\epsilon\,k^{\mathrm{d}})+
O\left(k^{-\infty}\right).
\end{eqnarray*}
Hence, $\lim_{k\rightarrow +\infty}(\pi/k)^{\mathrm{d}}\,\dim\left(V_{k\lambda}^{(k)}\right)
=\mathrm{vol}\big(M_{<\lambda}\big)+O(\epsilon)$, for any $\epsilon>0$.
Letting $\epsilon\rightarrow 0^+$, we get the statement.

\hfill Q.E.D.

\begin{rem}
By the same argument, the weak limit of
$(\pi/k)^{\mathrm{d}}\,\mathfrak{T}_k\left(k\lambda,\cdot\right)$
as $k\rightarrow +\infty$
is the characteristic function of $M_{<\lambda}$.
\end{rem}

\section{Appendix}
Although not strictly necessary, let us briefly pause give a direct derivation of the
asymptotic expansion which, in view of (\ref{eqn:2-nd-computation-key-integral})
below, is equivalent to the scaling limit discussed in the introduction (cfr  \cite{sz}).
Let us now suppose $\lambda>\varsigma_T(x)$. Before we proceed the following remark
is in order.

Let $\nu>0$ be such that $\lambda>\varsigma_T(x)+3\,\nu$, and choose
$g\in \mathcal{C}^\infty_0(\mathbb{R})$ such that
$g\ge 0$, $g(b)=1$ if $b\in \big(\varsigma_T(x)-\nu,\varsigma_T(x)+\nu\big)$,
$g(b)=0$ if $b\not\in \big(\varsigma_T(x)-2\nu,\varsigma_T(x)+2\nu\big)$.
Multiplying the integrand in (\ref{eqn:1-st-integral-form-F-stainary-phase})
by the identity $1=g(b)+\big(1-g(b)\big)$ the integral splits as the sum of two terms.
In the second of these, the integrand is supported where $\big|b-\varsigma_T(x)\big|
\ge\nu$; by Proposition \ref{prop:rapid-decay-case},
as $k\rightarrow +\infty$ this is rapidly decreasing.

Therefore we need only worry about the first summand. We may thus
assume that the integrand is compactly supported in $b$, tacitly
absorb the cut-off in the amplitude, and rewrite (\ref{eqn:1-st-integral-form-F-stainary-phase})
as

\begin{eqnarray}
\label{eqn:1-st-integral-form-F-stationary-compact-supp}
\lefteqn{\int _{-\infty}^{+\infty}G(k\lambda-\eta)\, d\mu_{\mathcal{T}^{(k)}}(\eta)\sim \frac{k^{\mathrm{d}+3}}{(2\pi)^{2\mathrm{d}+2}}}\\
&&\cdot
\int_{\mathbb{C}^\mathrm{d}}\int_{S^{2\mathrm{d}}}e^{-i\sqrt{k}\,\mathbf{u}\cdot \omega_1}
\left[\int_1^{D}\int _{-\varepsilon}^\varepsilon\int_{-\pi}^\pi
\int_{1/2}^\lambda \int_{1/C}^C\int_{-\epsilon}^\epsilon
e^{ik\Phi}\,J_k\, d\xi\right]\,d\mathbf{u}\,d\omega,\nonumber
\end{eqnarray}
where $d\xi=dt\,d\theta\,d\vartheta\,db\,dr\,d\tau$.
Integration in $d\mathbf{u}$ is over a ball of radius
$O\left(k^{1/6}\right)$ in $\mathbb{C}^\mathrm{d}$.

We are now in a position to apply the stationary phase Lemma to determine
the asymptotics of the inner integral in (\ref{eqn:1-st-integral-form-F-stationary-compact-supp}),
by viewing $\mathbf{u}\in \mathbb{C}^\mathrm{d}$ and $\omega\in S$ as parameters.
A straightforward computation then leads to the following:

\begin{lem}
\label{lem:stationary-point}
For every $\omega\in S=S_\delta$,
$\Phi=\Phi(t,\theta,\vartheta,b,r,\tau)$ has the unique stationary point
$R_0=(t_0,\theta_0,\vartheta_0,b_0,r_0,\tau_0)=\big(1,0,0,q(x,\omega)/\omega_0,1/\omega_0,0\big)$.
At $R_0$, the Hessian of $\Phi$ has determinant $\det\left(\Phi''(R_0)\right)=-\omega_0^2$.
\end{lem}

Since $\omega_0>\delta>0$ on $S$, the stationary point is always non degenerate.
%
Furthermore, $\Phi(R_0)=0$ and $\det\big(k\,\Phi''(R_0)/2\pi i\big)^{1/2}=\big(k/2\pi \big)^3\,\omega_0$.
Let us make the assumption $\chi(0)=1$.
In view
of Remark \ref{rem:leading-order-term}, the stationary phase Lemma
implies
that the inner integral in (\ref{eqn:1-st-integral-form-F-stationary-compact-supp})
is given by
\begin{eqnarray}
\frac{8}{\omega_0}\cdot \left(\frac{k}{\pi}\right)^{\mathrm{d}-3}
\,e^{\psi_2(\omega_0 \mathbf{u},\mathbf{v})}\,\left(1+\sum_{j= 1}^Nk^{-j/2}\,c_j\right)+R_N,\nonumber
\end{eqnarray}
where
$\big\|R_N(\mathbf{u},\omega,\mathbf{v})\big\|_{C^j}
\le C_N\,e^{-(1-\upsilon)\|(\omega_0 \mathbf{u})-\mathbf{v}\|^2/2}\,k^{-(N+1)/2}$
(here $0<\upsilon\ll 1$)
on the domain of integration in $(\mathbf{u},\omega_1)$.

Writing $\mu=\sqrt{k}$ we have
\begin{eqnarray}
\label{eqn:1-st-integral-form-F-term-by-term}
\lefteqn{\int _{-\infty}^{+\infty}G(k\lambda-\eta)\, d\mu_{\mathcal{T}^{(k)}}(\eta)\sim
\frac{1}{\pi^{\mathrm{d}}}\cdot \frac{\mu^{4\mathrm{d}}}{(2\pi)^{2\mathrm{d}-1}}\,}\\
&&\cdot
\left\{\int_{\mathbb{C}^\mathrm{d}}\int_{S^{2\mathrm{d}}}e^{-i\mu\,\mathbf{u}\cdot \omega_1}
\left[\left(\frac{e^{\psi_2(\omega_0 \,\mathbf{u},\mathbf{v})}}{\omega_0}\right)
\,
\left(1+\sum_{j= 1}^N\mu^{-j}\,c_j\right)+
R_N(\mathbf{u},\omega,\mathbf{v})\right]\,d\mathbf{u}\,d\omega\right\}.\nonumber
\end{eqnarray}
We have a unique stationary point at $\mathbf{u}=\omega_1=\mathbf{0}$.
At this point the Hessian has determinant one, and therefore
\begin{eqnarray}
\label{eqn:1st-comput-final-exp}
\lefteqn{\int _{-\infty}^{+\infty}G(k\lambda-\eta)\, d\mu_{\mathcal{T}^{(k)}}(\eta)}\\
&\sim& \frac{1}{\pi^{\mathrm{d}}}\cdot \frac{\mu^{4\mathrm{d}}}{(2\pi)^{2\mathrm{d}-1}}\,\,
\left(\frac{2\pi}{\mu}\right)^{2\mathrm{d}}\,e^{-\|\mathbf{v}\|^2/2}
\,\left(1+\sum_{j\ge 1}\mu^{-j}\,d_j\right)\nonumber\\
&=&2\pi\,\left(\frac k\pi\right)^{d}\,e^{-\|\mathbf{v}\|^2/2}
\,\left(1+\sum_{j\ge 1}\mu^{-j}\,d_j\right).\nonumber
\end{eqnarray}
The previous expression also gives the correct bound on the remainder.


\begin{thebibliography}{Dillo99}


\bibitem[BPU]{bpu} D. Borthwick, T. Paul, A. Uribe, {\em
Semiclassical spectral estimates for Toeplitz operators},
Ann. Inst. Fourier (Grenoble) \textbf{48} (1998), no. 4, 1189--1229


\bibitem[BSZ]{bsz} P. Bleher, B. Shiffman, S. Zelditch, {\em
Universality and scaling of correlations between zeros on complex
manifolds}, Invent. Math. \textbf{142} (2000), 351--395


\bibitem[BG]{bg} L. Boutet de Monvel, V. Guillemin, {\em The spectral theory of Toeplitz operators}, Annals of Mathematics Studies, \textbf{99} (1981), Princeton University Press, Princeton, NJ; University of Tokyo Press, Tokyo


\bibitem[BS]{bs} L. Boutet de Monvel, J. Sj\"ostrand,
{\em Sur la singularit\'e des noyaux de Bergman et de Szeg\"o},
Ast\'erisque \textbf{34-35} (1976), 123--164



\bibitem[C]{c} M. Christ,
{\it Slow off-diagonal decay for Szeg\"o kernels associated to smooth Hermitian line bundles},
Harmonic analysis at Mount Holyoke (South Hadley, MA, 2001),
77--89, Contemp. Math. \textbf{320}, Amer. Math. Soc., Providence,
RI, 2003




\bibitem[GS]{gs} A. Grigis, J. Sj\"ostrand, {\em Microlocal analysis for differential operators. An introduction},
London Mathematical Society Lecture Note Series,  \textbf{196} (1994),
Cambridge University Press, Cambridge

\bibitem[H]{h} L. H\"{o}rmander, {\em
The spectral function of an elliptic operator}, Acta Math.  \textbf{121} (1968), 193--218


\bibitem[P]{p} R. Paoletti, {\em
On the Weyl law of a Toeplitz operator}, preprint

\bibitem[SZ]{sz} B. Shiffman, S. Zelditch, {\em Asymptotics of almost
holomorphic sections of ample line bundles on symplectic
manifolds}, J. Reine Angew. Math. {\bf 544} (2002), 181--222

\bibitem[T]{t} G. Tian,
{\em On a set of polarized K\"{a}hler metrics on algebraic manifolds}, J. Differential Geom. \textbf{32} (1990), no. 1, 99--130

\bibitem[Z]{z} S. Zelditch, {\em Szeg\"o kernels and a theorem of Tian},
Int. Math. Res. Not. {\bf 6} (1998), 317--331

\end{thebibliography}
\end{document}